\documentclass[12pt]{article}
\usepackage{filecontents}

\usepackage{amssymb}
\usepackage{amsthm}
\usepackage{amssymb}
\usepackage{rotating}
\usepackage{caption}
\usepackage{array}
\usepackage{color}
\usepackage{multirow}
\usepackage{mathrsfs}
\usepackage{subfigure}
\usepackage{graphicx,epstopdf,float}
\usepackage{subfigure}
\usepackage{amsmath,bm}
\usepackage{placeins}
\usepackage{enumerate}
\usepackage[margin=0.5in]{geometry}
\usepackage[title]{appendix}
\usepackage{enumitem}
\usepackage{varioref}
\usepackage{indentfirst}
\usepackage{setspace}
\usepackage{amsfonts}
\usepackage{booktabs,tabularx}
\usepackage{lscape}
\usepackage{xcolor}
\usepackage{lineno}
\usepackage{pgfplots}
\usepackage{dcolumn}  
\usepackage{array} 
\usepackage{algorithm}%
\usepackage{algorithmicx}%
\usepackage[authoryear]{natbib}
\usepackage{hyperref}
\usepackage{secdot}
\usepackage{apalike}
\usepackage{graphicx}%
\usepackage{multirow}%
\usepackage{amsmath,amssymb,amsfonts}%
\usepackage{amsthm}%
\usepackage{mathrsfs}%
\usepackage[title]{appendix}%
\usepackage{xcolor}%
\usepackage{textcomp}%
\usepackage{manyfoot}%
\usepackage{booktabs}%
\usepackage{algorithm}%
\usepackage{algorithmicx}%
\usepackage{algpseudocode}%
\usepackage{listings}%
\usepackage{relsize}
\usepackage{caption}

\usepackage[utf8]{inputenc}
\usepackage{array}
\usepackage{makecell}

\begin{document}
	\newcommand{\bea}{\begin{eqnarray}}
		\newcommand{\eea}{\end{eqnarray}}
	\newcommand{\nn}{\nonumber}
	\newcommand{\bee}{\begin{eqnarray*}}
		\newcommand{\eee}{\end{eqnarray*}}
	\newcommand{\lb}{\label}
	\newcommand{\nii}{\noindent}
	\newcommand{\ii}{\indent}
	\newtheorem{theorem}{Theorem}[section]
	\newtheorem{longtable}{longtable}
	\newtheorem{example}{Example}[section]
	\newtheorem{corollary}{Corollary}[section]
	\newtheorem{definition}{Definition}[section]
	\newtheorem{lemma}{Lemma}[section]
	\newtheorem{dataset}{Dataset}[section]
	\newtheorem{Nomenclature}{Nomenclature}
	\newtheorem{remark}{Remark}[section]
	\newtheorem{proposition}{Proposition}[section]
	\numberwithin{equation}{section}
	\renewcommand{\theequation}{\thesection.\arabic{equation}}
	\bibpunct[, ]{(}{)}{;}{a}{,}{,}
	\renewcommand\bibfont{\fontsize{10}{12}\selectfont}
	\setlength{\bibsep}{0.0pt}
	\title{\bf Block Adaptive Progressive Type-II Censored Sampling for the Inverted Exponentiated Pareto Distribution: Parameter Inference and Reliability Assessment}
	\author{Rajendranath Mondal$^a$\thanks{Email address: rajendranathmondal9@gmail.com, rajendranath\_m@ma.iitr.ac.in}, Aditi Kar Gangopadhyay$^a$\thanks{Email address (corresponding author): aditi.gangopadhyay@ma.iitr.ac.in}, Raju Bhakta$^a$\thanks{Email address: bhakta.r93@gmail.com, raju.cs@sric.iitr.ac.in}~~and Kousik Maiti$^b$\thanks{Email address: kousikulu@gmail.com}
	\\{\it\small$^a$Department of Mathematics, Indian Institute of Technology Roorkee, Roorkee-247667, Uttarakhand, India.}
	\\{\it\small$^b$School of Applied Science $\&$ Humanities, Haldia Institute of Technology, Haldia-721657, West Bengal, India.}}
	\date{}
	\maketitle 
    \begin{abstract}   
    This article explores the estimation of unknown parameters and reliability characteristics under the assumption that the lifetimes of the testing units follow an Inverted Exponentiated Pareto (IEP) distribution. Here, both point and interval estimates are calculated by employing the classical maximum likelihood and a pivotal estimation methods. Also, existence and uniqueness of the maximum likelihood estimates are verified. Further,  asymptotic confidence intervals are derived by using the asymptotic normality property of the maximum likelihood estimator.  Moreover, generalized confidence intervals are obtained by utilizing the pivotal quantities. Additionally, some mathematical developments of the IEP distribution are discussed based on the concept of order statistics. Furthermore, all the estimations are performed on the basis of the block censoring procedure, where an adaptive progressive Type-II censoring is employed to every block. In this regard, the performances of two estimation methods, namely maximum likelihood estimation and pivotal estimation, is evaluated and compared through a simulation study. Finally, a real data is illustrated to demonstate the flexibility of the proposed IEP model.
    \end{abstract}	

    \section{Introduction}
    To evaluate the reliability characteristics (RCs) of the products, life testing experiments are conducted. In practice, censoring schemes are employed when the experiment is costly enough and consumes a large amount of time. In broad sense, censoring schemes are predetermined plans which are implemented to the life tests to create censored or incomplete lifetime data. Type-I and Type-II censorings are very well-known traditional censoring schemes and they are used by numerous authors. However, due to some drawbacks of Type-I and Type-II censorings, many other types of censoring schemes were devised. Among them, progressive censorings are significantly popular. Progressive censorings allow removals of working products during the experiment. Progressive Type-I and progressive Type-II are the simplest progressive censoring schemes. For detailed study of progressive Type-I and Type-II censorings, we refer \cite{balakrishnan2014art}. \cite{cramer2010adaptive} introduced an adaption on the progressive Type-II censoring and they named it as `adaptive  progressive Type-II censoring'. The authors established that maximum likelihood estimators obtained by using adaptive  progressive Type-II censoring coincide with the same obtained by using progressive Type-II censoring. Also, some inferential results on one and two-parameter exponential distribution were presented in this article. \cite{ng2009statistical} introduced a particular kind of adaption on the progressive Type-II censorings, which can be considered as a special case of the adaptive progressive Type-II censoring introduced by \cite{cramer2010adaptive}. In this article, point and interval estimations for the failure rate parameter of the exponential distribution were performed based on the proposed censoring scheme and the computaions for expected total test time were performed. 
    
    \cite{sobhi2016estimation} studied the estimation of the unknown parameters and the reliability and hazard functions of the two-parameter exponentiated Weibull distribution under the adaptive progressive Type-II censoring. Here, both the maximum likelihood and Bayesian estimations are performed. Markov chain Monte Carlo method and importance sampling methods are used to compute the Bayesian estimates and credible intervals for the unknown quantities. Finally, the comparisons between the proposed estimation methods are discussed.
        \cite{sewailem2019inference} estimated the unknown parameters and the survival function of the log-logistic distribution by employing both the maximum likelihood and Bayesian inference method based on the adaptive  progressive Type-II censored samples. The authors compared the proposed methods by performing a Monte-Carlo study.
        Further, \cite{nassar2023estimation} used an improved adaptive  progressive Type-II censoring to delve into the estimations of the unknown parameters, reliability and hazard functions  of the Weibull distribution. Here, the authors performed the maximum likelihood and maximum product of spacing estimation methods as classical approaches. Bayesian estimation was also done by applying the Markov chain Monte Carlo method. Finally, the proposed estimation methods are compared through a simulation study.  
     \cite{dutta2024bayesian} investigated inferences regarding the unknown parameters, as well as the reliability and hazard functions, within the framework of the adaptive progressive Type-II censoring scheme. Maximum likelihood estimates and the maximum product spacing estimates are derived numerically since they could not be obtained explicitly. Bayesian estimations for the same unknown quantities was executed and Markov chain Monte Carlo method was used since the Bayesian estimates could not be derived explicitly.

      \cite{abu2021statistical} introduced adaptive-general progressive Type-II censoring as a generalization of the adaptive progressive Type-II censoring. They investigated the inferences on the two-parameter Gompertz distribution by executing both the maximum likelihood and Bayesian methods of estimation. They performed a vast simulation study and analyzed a real data to compare the performances of these methods. There are other numerous authors who worked on statistical inferences on lifetime distributions using the adaptive  progressive Type-II censoring scheme. A significant drawback of all the censoring schemes mentioned previously in this section can be observed when the equipments required to test the lifetimes of all the units are not readily available. Moreover, sometimes it is not possible to observe all the failue times by only one observer. These situations can be partly overcome when the block censoring mechanisms are implemented. In these mechanisms, the collection of testing units are randomly partitioned into a prefixed finite number of groups (also known as facilities) and each group contains a prefixed number of units. Also, number of units in any two different groups may be equal or may be unequal. There are various kinds of block censoring mechanisms available in the literature.  \cite{johnson1964theory} introduced a censoring scheme, named first-failure censoring scheme (FFCS), in which units are divided into several groups and each group contains equal number of units. These groups can be tested simultaneously or consecutively. The test is terminated after the occurance of first failure in each test group. \cite{wu2009estimation} combined the progressive Type-II censoring scheme with FFCS to introduce the scheme named, progressive first-failure censoring scheme (PFFCS). Under this new censoring scheme, maximum likelihood estimation, exact and approximate confidence intervals and exact confidence region for the unknown parameters of the Weibull distribution are discussed in this article. Simulation studies were presented here to examine the performances of the proposed estimation methods. The expected total required time of the life test under PFFCS was also derived.

      \cite{ahmadi2018block} introduced a block censoring scheme in which Type-II censoring was executed in every facility. The authors implemented this censoring mechanism to estimate the parameters of the two-parameter exponential distribution. Also, the optimal block censoring was derived by a Monte Carlo simulation. \cite{kumari2023reliability} introduced the block progressive Type-II censoring scheme (BPCS) which is obtained by employing the progressive Type-II censoring instead of the traditional Type-II censoring in each facility of the block censoring setup mentioned in \cite{ahmadi2018block}. They estimated the RCs of a distribution introduced by \cite{chen2000new} under the block progressive Type-II censoring. They also estimated the differences of different test facilities of BPCS. The estimation methods used in this article are  maximum likelihood estimation and hierarchical Bayesian estimation. Approximate confidence intervals are constructed under the asymptotic theory and the delta method. Bayesian estimates are calculated through a hybrid Metropolis-Hastings sampling in hierarchical framework. Here, the proposed methods are compared through a simulation study. 
    
    In our article, we study various estimations on the two-parameter IEP distribution based on a block adaptive progressive Type-II censoring schemes (BAPCS). BAPCS is simply formed by replacing the progressive Type-II censoring by the adaptive progressive Type-II censoring in the BPCS setup. As the adaptive progressive censoring, we have used the model introduced by \cite{ng2009statistical}. Note that, \cite{ghitany2014likelihood} established existence and uniqueness of the maximum likelihood estimates (MLE) for a family of two-parameter inverted exponentiated distributions having the cumulative distribution function (CDF) of the form
    
    \begin{equation} \label{eq1}
    	F(t; \alpha, \beta) = 1-  \left(1 - \exp(-\beta Q(1/t))\right)^\alpha; ~~ t > 0, ~~ \alpha, \beta >0,
    \end{equation}
    where $Q$ is an increasing function satisfying $Q(0) = 0$ and $Q(+\infty) = +\infty$, and $\alpha$ and $\beta$ are both shape parameters of this family. The authors discussed the MLE methods for the unknown parameters $\alpha$ and $\beta$ of this family of distributions given in (\ref{eq1}), considering the complete sample as well as progressive Type-I and Type-II censored samples.
    
    Now, from (\ref{eq1}), the IEP distribution is derived by using the substitution $Q(1/t) = \log(1 + 1/t)$ to (\ref{eq1}). Consequently, the probability density function (PDF) of the IEP distribution is as follows: 
    \begin{equation}\label{eq2} \displaystyle{
    f(t; \alpha, \beta) = \alpha\beta t^{\beta - 1}(1 + t)^{-(\beta + 1)}\left[1 - \left(\frac{t}{1+t} \right)^\beta \right]^{\alpha - 1}; ~~ t > 0, ~~ \alpha, \beta >0. }
    \end{equation}

    Moreover, the reliability and hazard functions of the IEP distribution are 
    
    \begin{equation}
    R(t; \alpha, \beta) = \left[1 - \left(\frac{t}{1+t} \right)^\beta \right]^{\alpha} ; ~~ t > 0, ~~ \alpha, \beta >0
    \end{equation}

    and
    
    \begin{equation}
    \displaystyle{H(t; \alpha, \beta) = \frac{\alpha\beta t^{\beta - 1}(1 + t)^{-(\beta + 1)}}{1 - \left(\frac{t}{1+t} \right)^\beta} ; ~~ t > 0, ~~ \alpha, \beta >0},
    \end{equation}
    respectively. Furthermore, the median time to failure (MTF) of the IEP distribution is as follows: 
    
    \begin{equation}
    \mu(\alpha, \beta) = \left\{(1 - 2^{-1/\alpha})^{-1/\beta} - 1 \right\}^{-1}; ~~ \alpha, \beta >0.
    \end{equation}

  \cite{maurya2019inference} estimated the parameters of the IEP distribution by utilizing both the MLE method and Bayesian estimation method under the progressive Type-II censoring. Here, MLEs are derived by using the expectation-maximization algorithm and the Bayesian estimation was carried out by the importance sampling method based on both symmetric and asymmetric loss functions. Further, highest posterior density intervals as well as one and two sample Bayesian prediction were discussed. Moreover, the authors presented a simulation study to examine the performances of the proposed estimation and prediction methods. Also, the optimal censoring plans are devised based on two information measure criteria. Lastly, a real data analysis was done for illustration purpose. \cite{kumari2023rel} done a reliability estimation in a multicomponent system based on Type-II censored data, where the stress and strength variables follow the inverted exponentiated family of distributions introduced by \cite{ghitany2014likelihood} with the common scale parameter $\beta$. The authors performed the maximum likelihood estimation and approximate confidence interval estimation of the stress-strength reliabity function. Then, some pivotal quantities were used to obtain the pivotal estimate of this function. Also, the authors discussed the particular situation when all the parameters of the stress-strength components are unknown and various estimations for the reliability were suggested for this case. Finally, all the estimation procedures are compared via a simulation study.
  
  In our BAPCS procedure, considering the natural similarities among the test facilities in practical uses, we assume that the parameter $\beta$ is shared by all the facilities. Here, we characterize the differences in different test facilities (DDTF) by the distinct parameters $\alpha_i$ ($i=1,2,\cdots,k$), at the $i$-th test facility, where $k$ is the number of facilities in our BAPCS setup.
  In other words, we assume that the lifetimes of the units in the $i$-th facility of BAPCS setup independently follow IEP($\alpha_i, \beta$) distribution ($i = 1,2,...,k$). If the differences in different test facilities are not negligible, the population parameter $\alpha$ in  (\ref{eq2}) can be calculated as $\alpha = \sum_{i=1}^{k}\omega_i\alpha_i/\sum_{i=1}^{k}\omega_i$, where $\omega_i$'s are some weight coefficients.
  
  The rest of the paper is organized as follows. The Section \ref{sec 2} consists of the description of the MLE of the unknown parameters as well as that of the reliability function $R(t; \alpha, \beta)$,  hazard function $H(t; \alpha, \beta)$, and median time to failure $\mu(\alpha, \beta)$ of the IEP distribution. In Section \ref{sec 3}, the constructions of the asymptotic confidence intervals of $\alpha_i$ ($i = 1, 2, \cdots, k$), $\alpha$, $\beta$ are carried out. Also, here, the asymptotic confidence intervals of $R(t; \alpha, \beta)$, $H(t; \alpha, \beta)$, and $\mu(\alpha, \beta)$ are described. The pivotal estimation for estimating the same quantities are discussed thoroughly in Section \ref{sec 4}. In Section \ref{section5}, some of the mathematical developments of the IEP distribution are described based on the concept of order statistics. An extensive simulation study to compare the performances of the MLE method and the pivotal estimation method is presented in Section \ref{sec 6}. To illustrate the flexibility of the IEP model, a real data is analyzed in Section \ref{sec 7}. Finally, we present several important conclusions from our work in Section \ref{sec 8}.

  \section{Maximum Likelihood Estimation}\label{sec 2}
   
   Suppose, in a block censoring setup, $n$ number of life testing units are divided into $k$ number of groups (or facilities) and the $i$-th group contains $n_i$ number of units ($i = 1, 2, \dots,k$) such that $\sum_{i=1}^{k} n_i = n$. Let each facility is tested under the adaptive progressive model introduced by \cite{ng2009statistical}. Also, assume that the integers $m_i$ ($i = 1, 2, \dots ,k$) denote the number of observable failures in the $i$-th facility. If $\mathcal{T}_{ij}$ be the $j$-th failure time for the $i$-th facility ($j = 1,2, \dots, m_i$ and $i = 1, 2, \dots, k$), $T_i$ be the threshold time associated with $i$-th facility, and $J_i$ be the number of failures before the threshold time $T_i$ in the $i$-th facility, then as per the likelihood function for the adptive progressive Type-II censoring derived in \cite{ng2009statistical}, we can write the likelihood function for the $i$-th facility as 
   \begin{equation}
   	\mathcal{L}_i(\mathcal{T}_i; \alpha_i, \beta) = d_{J_i}\left(\prod_{j=1}^{m_i}f(\mathcal{T}_{ij}; \alpha_i, \beta)\right)\left(\prod_{j=1}^{J_i}(1 - F(\mathcal{T}_{ij}; \alpha_i, \beta))^{R_{ij}}\right)\left(1 - F(\mathcal{T}_{im_i}; \alpha_i, \beta)\right)^{n_i - m_i - \sum_{j=1}^{J_i}R_{ij}}.
   	\label{lkl}
   \end{equation}

  In (\ref{lkl}), the constants $d_{J_i}$ is defined as $d_{J_i} = \prod_{j=1}^{m_i}\left[ n_i - i + 1 - \sum_{k=1}^{\min\{j-1,J_i\}} R_{ik}\right]; ~~ i = 1, 2, \dots, k$ and the functions $F$ and $f$ are defined in (\ref{eq1}) and (\ref{eq2}), respectively. Also, $R_{ij}$'s in (\ref{lkl}), are the number of working units which are withdrawn at the $j$-th failure time point in the $i$-th facility. The random vector  $\mathcal{T}_{i} = (\mathcal{T}_{i1}, \mathcal{T}_{i2},\cdots,\mathcal{T}_{im_i})$ denotes the adaptive progressive Type-II censored failure time sample from the IEP($\alpha_i, \beta$) distribution situated in the $i$-th facility for $i = 1,2,\cdots, k$.  Although, the values of $R_{ij}$'s are preassumed, removing of working units during the experiment is done as per the censoring model described in \cite{ng2009statistical}. Consequently, the likelihood function for the entire block adaptive progressive Type-II setup is given as

   \begin{equation}
   	\begin{split}
   \mathcal{L}(\mathcal{T}; \alpha_1, \alpha_2, \dots, \alpha_k, \beta) =	\prod_{i=1}^{k} d_{J_i}\left(\prod_{j=1}^{m_i}f(\mathcal{T}_{ij}; \alpha, \beta)\right)\left(\prod_{j=1}^{J_i}(1 - F(\mathcal{T}_{ij}; \alpha_i, \beta))^{R_{ij}}\right)\\ \times \left(1 - F(\mathcal{T}_{im_i}; \alpha_i, \beta)\right)^{n_i - m_i - \sum_{j=1}^{J_i}R_{ij}},
   \end{split}
  \end{equation} 

  where $\mathcal{T} = (\mathcal{T}_{1},\mathcal{T}_{2},\cdots, \mathcal{T}_{k})$ denotes the combined BAPCS samples from all the facilities.
  Hence, the log-likelihood function, after using the functional forms of $f$ and $F$, and omitting the terms $d_{J_i}$'s,  is as follows: 
  
  \begin{equation}
  \begin{split} \label{llh}
  l(\mathcal{T}; \alpha, \beta) = \sum_{i=1}^k \sum_{j=1}^{m_i}\left\{ \log(\alpha_i) + log(\beta) + (\beta - 1)\log(\mathcal{T}_{ij}) - (\beta + 1)\log(1 + \mathcal{T}_{ij})  \right\} \\ + \sum_{i=1}^k \sum_{j=1}^{m_i} (\alpha_i - 1)  \log\left\{ 1 - \left( \frac{\mathcal{T}_{ij}}{ 1 + \mathcal{T}_{ij}} \right)^\beta \right\} + \sum_{i=1}^k \sum_{j=1}^{J_i} \alpha_i R_{ij} \log\left\{ 1 - \left( \frac{\mathcal{T}_{ij}}{ 1 + \mathcal{T}_{ij}} \right)^\beta \right\}\\ + \sum_{i=1}^{k} \alpha_i \left(n_i - m_i - \sum_{j=1}^{J_i}R_{ij}\right)\log\left\{ 1 - \left( \frac{\mathcal{T}_{im_i}}{ 1 + \mathcal{T}_{im_i}} \right)^\beta \right\}.
  \end{split}
  \end{equation}

  Note that, with a proper application of the theorem and lemmas mentioned in \cite{ghitany2014likelihood}, it can be shown that the log-likelihood function given in (\ref{llh}) has unique global maximum at \\
  
   $\displaystyle \alpha_i = - \frac{m_i}{ \sum_{j=1}^{m_i} \log\left( 1 - t_{ij}^{\hat{\beta}} \right) + \sum_{j=1}^{J_i}R_{ij}\log\left( 1 - t_{ij}^{\hat{\beta}} \right) + R_{im_i} \log\left( 1 - t_{im_i}^{\hat{\beta}} \right)} = \hat{\alpha}_i \mbox{ (say)}$
   
   and \\
   
   $\beta = \hat{\beta}$, where $\hat{\beta}$ is the unique solution of the equation 
   
   \begin{equation}
   	\label{bteq}
   	\begin{split}
   \frac{\sum_{i=1}^{k}m_i}{\beta} + \sum_{i=1}^k \sum_{j=1}^{m_i} \log\left( t_{ij} \right) + \sum_{i=1}^k \sum_{j=1}^{m_i}\frac{  t_{ij}^\beta \log(t_{ij})  }{1 - t_{ij}^\beta} ~~ ~~ ~~ ~~ ~~ ~~ ~~ ~~ ~~ ~~ ~~ ~~ ~~ ~~ ~~ ~~ ~~ ~~ ~~ ~~ ~~ ~~ ~~ ~~ ~~ ~~ ~~ ~~ ~~ ~~ ~~ ~~ ~~ ~~ ~~ ~~ ~~ ~~ ~~ \\ + \sum_{i=1}^{k}m_i \frac { \sum_{j=1}^{m_i}t_{ij}^\beta \log\left( t_{ij}\right)/\left(1 -t_{ij}^\beta \right)+\sum_{j=1}^{J_i}R_{ij}\left( t_{ij} \right)^\beta \log\left( t_{ij} \right)/\left(1 - t_{ij}^\beta \right) + R_{im_i}t_{im_i}^\beta \log\left( t_{im_i} \right)/(1-t_{im_i}^\beta)}{\sum_{j=1}^{m_i} \log\left(1 - t_{ij}^\beta \right) + \sum_{j=1}^{J_i}R_{ij}\log\left( 1 - t_{ij}^\beta \right) + R_{im_i} \log\left( 1 - t_{im_i}^\beta \right)} = 0,
   \end{split}  
   \end{equation}

   for the variable $\beta$, provided $t_{ij} = \frac{\mathcal{T}_{ij}}{1 + \mathcal{T}_{ij}}$ and $R_{im_i} = \left(n_i - m_i - \sum_{j=1}^{J_i}R_{ij}\right)$ ($j = 1,2, \dots,m_i$ and $i = 1, 2, \dots, k$). Therefore, from above justification, we obtain unique MLEs of $\alpha$ and $\beta$.
   Consequently, the MLEs of $\alpha_i$ and $\beta$ are $\hat{\alpha}_i$ and $\hat{\beta}$, respectively ($i = 1, 2, \dots, k$). It should be noted that, (\ref{bteq}) cannot be solved in closed form expressions and so, we have to solve it numerically.  Further, due to the invariance property of MLE, the MLE of any parametric function $\phi(\alpha, \beta)$ is $\phi(\hat{\alpha}, \hat{\beta})$. Therefore, the MLEs of the reliability function, hazard function, and MTF are obtained as

   \begin{eqnarray}
   \hat{R}(t) &=& \left[1 - \left(\frac{t}{1+t} \right)^{\hat{\beta}} \right]^{\hat{\alpha}},\\
   \hat{H}(t) &=& \frac{\hat{\alpha}\hat{\beta} t^{\hat{\beta} - 1}(1 + t)^{-(\hat{\beta} + 1)}}{1 - \left(\frac{t}{1+t} \right)^{\hat{\beta}}} \\
   \mbox{and}\\ \hat{\mu} &=& \mu(\hat{\alpha}, \hat{\beta}) = \left\{(1 - 2^{-1/\hat{\alpha}})^{-1/\hat{\beta}} - 1 \right\}^{-1},
   \end{eqnarray}

   respectively. Here, $\hat{\alpha} = \sum_{i=1}^{m}\hat{\omega_i}\hat{\alpha}_i/\sum_{i=1}^{m}\hat{\omega}_i$ is an weighted estimate of $\alpha$, where $\hat{\omega}_i = 1/\hat{V}_{\alpha_i}$ ($\hat{V}_{\alpha_i}$ is the observed variance of $\alpha_i$ mentioned in Section \ref{sec 3}) and $\hat{\alpha}$ is considered as the MLE of the parameter $\alpha$.
   
   \section{Asymptotic Confidence Intervals} \label{sec 3}
   
  In this section, the asymptotic interval estimations for the unknown parameters and RCs of the IEP distribution are discussed. Due to the asymptotic normality property of MLE (for detailed study for this property, we refer \cite{vander1990accuracy}), ($\hat{\beta}, \hat{\alpha}_1, \hat{\alpha}_2, \dots, \hat{\alpha}_k$) jointly follows $(k+1)$-variate normal distribution with mean vector ($\beta, \alpha_1, \alpha_2, \dots, \alpha_k$) and variance-covariance matrix $I^{-1}$, where $I$ is the Fisher information matrix defined as
  
  \begin{eqnarray}
  	I = -\begin{bmatrix}
  		E\left(\frac{\partial^2l}{\partial\beta^2}\right) & E\left(\frac{\partial^2l}{\partial\beta\partial\alpha_1}\right) & E\left(\frac{\partial^2l}{\partial\beta\partial\alpha_2}\right) & \dots & E\left(\frac{\partial^2l}{\partial\beta\partial\alpha_k}\right)\\
  		E\left(\frac{\partial^2l}{\partial\alpha_1\partial\beta}\right) & E\left(\frac{\partial^2l}{\partial\alpha_1^2}\right) & E\left(\frac{\partial^2l}{\partial\alpha_1\partial\alpha_2}\right) & \dots & E\left(\frac{\partial^2l}{\partial\alpha_1\partial\alpha_k}\right)\\
  		\vdots & \vdots & \vdots & \ddots & \vdots\\
  		E\left(\frac{\partial^2l}{\partial\alpha_k\partial\beta}\right) & E\left(\frac{\partial^2l}{\partial\alpha_k\alpha_1}\right) & E\left(\frac{\partial^2l}{\partial\alpha_k\partial\alpha_2}\right) & \dots & E\left(\frac{\partial^2l}{\partial\alpha_k^2}\right)\\
  	\end{bmatrix},
  \end{eqnarray}
where $l$ is the log-likelihood function defined in (\ref{llh}) and $E$ is the expectation defined over the joint distribuition of the random failure times $\mathcal{T}_{ij}$'s ($j=1,2,\dots,m_i; i = 1,2,\dots,k$). Here, the expectations taken in the Fisher information matrix is not possible to be obtained in closed form expressions. The entries of the Fisher information matrix is estimated by removing the expectations and evaluating the double derivatives at $\alpha_i = \hat{\alpha}_i$ and $\beta = \hat{\beta}$. Then, the estimated Fisher information matrix is given by
  
  \begin{eqnarray}
  	\hat{I} =
  	-\begin{bmatrix}
  		\frac{\partial^2l}{\partial\beta^2}& \frac{\partial^2l}{\partial\beta\partial\alpha_1} & \frac{\partial^2l}{\partial\beta\partial\alpha_2} & \dots & \frac{\partial^2l}{\partial\beta\partial\alpha_k}\\
  		\frac{\partial^2l}{\partial\alpha_1\partial\beta} & \frac{\partial^2l}{\partial\alpha_1^2} & \frac{\partial^2l}{\partial\alpha_1\partial\alpha_2} & \dots & \frac{\partial^2l}{\partial\alpha_1\partial\alpha_k}\\
  		\vdots & \vdots & \vdots & \ddots & \vdots\\
  		\frac{\partial^2l}{\partial\alpha_k\partial\beta} & \frac{\partial^2l}{\partial\alpha_k\alpha_1} & \frac{\partial^2l}{\partial\alpha_k\partial\alpha_2} & \dots & \frac{\partial^2l}{\partial\alpha_k^2}\\
  	\end{bmatrix}_{(\beta = \hat{\beta},  \alpha_1 = \hat{\alpha}_1,  \alpha_2 = \hat{\alpha}_2,\dots,  \alpha_k = \hat{\alpha}_k)},
  \end{eqnarray}
which is used to derive the asymptotic confidence invervals of the unknown parameters as well as of the RCs. In fact, if $\hat{V} = \hat{I}^{-1}$ and the entries of $\hat{V}$ be denoted as $v_{ij}$ ($i,j = 1,2,\dots,k +1$), then, $v_{11}$ is considered as the estimated asymptotic variance of $\hat{\beta}$ and $v_{ii}$ is considered as the estimated asymptotic variance of $\alpha_{i-1}$ for $i = 2,3,\dots,k+1$. Consequently, the $(1-\gamma)\%$ confidence interval for $\beta$ is given as

 \begin{align} \nonumber \left(\hat{\beta} - z_{\gamma/2}\sqrt{v_{11}}, ~~  \hat{\beta} + z_{\gamma/2}\sqrt{v_{11}}\right)  
\end{align} 
and the same for $\alpha_{i-1}$ is as follows:
\begin{align} \nonumber  \left(\hat{\alpha}_{i-1} - z_{\gamma/2}\sqrt{v_{ii}}, ~~  \hat{\alpha}_{i-1} + z_{\gamma/2}\sqrt{v_{ii}}\right); ~~ (i = 2,3,\dots,k+1), 
\end{align}
   where $z_{\gamma/2}$ is defined as $P(Z>z_{\gamma/2}) = \gamma/2$, $Z$ follows standard normal distribution, i.e., $Z \sim N(0, 1)$.
   Further, to construct a confidence interval for any parametric function $\phi(\beta,\alpha_1,\alpha_2,\dots,\alpha_k)$, we should make use of the delta method (see \cite{beutner2024delta}). According to this method, $\phi(\hat{\beta},\hat{\alpha}_1,\hat{\alpha}_2,\dots,\hat{\alpha}_k)$ also has asymptotic normal distribution with mean $\phi(\beta,\alpha_1,\alpha_2,\dots,\alpha_k)$ and variance\\ $V_{\phi} = [\nabla \phi(\beta,\alpha_1, \alpha_2, \dots, \alpha_k)]^{T} I^{-1}[\nabla \phi(\beta,\alpha_1, \alpha_2, \dots, \alpha_k)]$, where $I^{-1}$ is the inverse of the Fisher information matrix and $\nabla$ is the gradient operator. Hence, the variance $V_{\phi}$ is estimated as\\
   $\hat{V}_{\phi} = [\nabla \phi(\hat{\beta},\hat{\alpha}_1, \hat{\alpha}_2, \dots, \hat{\alpha}_k)]^{T} \hat{V}[\nabla \phi(\hat{\beta},\hat{\alpha}_1, \hat{\alpha}_2, \dots, \hat{\alpha}_k)]$ and then, the asymptotic confidence interval of $\phi(\hat{\beta},\hat{\alpha}_1,\hat{\alpha}_2,\dots,\hat{\alpha}_k)$ is obtained as 
   
    \begin{align}\left(\phi(\hat{\beta},\hat{\alpha}_1, \hat{\alpha}_2, \dots, \hat{\alpha}_k) - z_{\gamma/2}\sqrt{\hat{V}_{\phi}}, ~~  \phi(\hat{\beta},\hat{\alpha}_1, \hat{\alpha}_2, \dots, \hat{\alpha}_k) + z_{\gamma/2}\sqrt{\hat{V}_{\phi}} \right)
    \tag{A}
    \label{int_ph}.
    \end{align}
 Now, the asymptotic confidence intervals of the RCs $R(t;\alpha,\beta)$, $H(t;\alpha,\beta)$ and $\mu(\alpha, \beta)$ should be made by replacing the function $\phi$ by these functions in the above expression (\ref{int_ph}) of the confidence interval of $\phi(\beta,\alpha_1,\alpha_2,\dots,\alpha_k)$. Here, it is provided that, $\alpha = \left(\sum_{i=1}^{k}\omega_i\alpha_i\right)/(\sum_{i=1}^{k}\omega_i)$ and the weight coefficients $\omega_i$'s are taken as $\omega_i = 1/\hat{V}_{\alpha_i}$, and $\hat{V}_{\alpha_i}$ is the entry $v_{(i+1)(i+1)}$ of the matrix $\hat{V}$, which is defined previously in this section ($i = 1,2,\dots,k$). Also, the asymptotic confidence interval for $\alpha$ is obtained by replacing the function $\phi$ by $\alpha$ (as a function of $\alpha_1,\alpha_2,\dots,\alpha_k$) to the above mentioned interval given in (\ref{int_ph}).

  \section{Pivotal Inference} \label{sec 4}
  In this section, inferences on the unknown parameters and RCs of the IEP distribution are studied by utilizing some certain pivotal quantities. Firstly, we construct the following transformation of random variables:
  
  $\mathcal{U}_{ij} = -\log\left(1 -F(\mathcal{T}_{ij};\alpha,\beta) \right) = -\alpha_{i}\log\left( 1 - \left( \frac{\mathcal{T}_{ij}}{1+\mathcal{T}_{ij}}\right)^\beta \right); ~~ (j = 1,2,\dots,m_i ~~\mbox{and}~~ i = 1, 2, \dots, k),$
  where $F$ is the CDF of the IEP distribution provided in (\ref{eq2}) and $\mathcal{T}_{ij}$'s are the failure times as defined in Section \ref{sec 2}. Now, it can be easily verified that ($\mathcal{U}_{i1}, \mathcal{U}_{i2}, \dots, \mathcal{U}_{im_i} $) forms an adaptive progressive Type-II censored sample of size $m_i$ drawn from the standard exponential distribution for every $i = 1,2,\dots,k$. Further, another transformation of random variables is considered as follows:
  \begin{eqnarray}{\nonumber}
  	\mathcal{V}_{i1} &=& \Gamma_{i1}\mathcal{U}_{i1},\\{\nonumber}
  	\mathcal{V}_{i2} &=& \Gamma_{i2}[\mathcal{U}_{i2} - \mathcal{U}_{i1}],\\{\nonumber}
  	\mathcal{V}_{i3} &=& \Gamma_{i3}[\mathcal{U}_{i3} - \mathcal{U}_{i2}],\\{\nonumber}
  	&\cdots&\\{\nonumber}
  	\mathcal{V}_{im_i} &=& \Gamma_{im_i}[\mathcal{U}_{im_i} - \mathcal{U}_{im_i-1}],
  \end{eqnarray}
where $\Gamma_{ij}$ is the number of existing units in the $i$-th facility of the BAPCS setup just before occurance of the $j$-th failure to that facility. Clearly, $\Gamma_{ij} = \left(n_i - j + 1 - \sum_{l=1}^{\min\{J_i,j-1\}} R_{il} \right).$ One can verify that, ($\mathcal{V}_{i1}, \mathcal{V}_{i2}, \dots, \mathcal{V}_{im_i} $) is a random sample of size $m_i$ drawn from the standard exponential distribution. Now, we define the random variables
  
  $\mathcal{W}_{ij} = \sum_{l=1}^{j}\mathcal{V}_{il} = - \alpha_i\sum_{l=1}^{j-1}(r_{il} + 1)\log\left[ 1 - \left(\frac{\mathcal{T}_{il}}{1 + \mathcal{T}_{il}}\right)^\beta \right] - \alpha_i\left( n_i - \sum_{l=1}^{j-1}(r_{il} +1) \right)\log\left[ 1 - \left(\frac{\mathcal{T}_{ij}}{1 + \mathcal{T}_{ij}}\right)^\beta \right]$, where the integers $r_{ij}$ should be defined in two diferent cases as follows:

  \textbf{Case I}: When $J_i < (m_i - 1)$, then
   \begin{equation}\nonumber
  	r_{ij} = \begin{cases}
  		R_{ij}, & \text{if $j\leq J_i$},\\
  		0, & \text{if $J_i < j \leq (m_i-1)$ }, i = 1, 2, \dots, k,
  	\end{cases}
\end{equation}

and

\begin{equation} \nonumber
	 r_{im_i} = \left(n_i - m_i - \sum_{l=1}^{J_i} R_{il}\right).
\end{equation}

\textbf{Case II}: When $J_i = (m_i - 1)$ or $m_i$, then
 
 \begin{equation}\nonumber
 r_{ij} = R_{ij}, ~~\mbox{for}~~ j = 1,2, \dots,m_i~~ \mbox{and} ~~i = 1, 2, \dots, k.
 \end{equation}

 Again, the ratios $\mathcal{Q}_{ij} = \mathcal{W}_{ij}/\mathcal{W}_{im_i};~~(j = 1,2,\dots,m_i - 1)$ form an order statistic of size $(m_i - 1)$ from the standard uniform distribution. Moreover, $\mathcal{Q}_{ij}$'s are independent of $$\mathcal{W}_{im_i} = - \alpha_i\sum_{l=1}^{m_i-1}(r_{il} + 1)\log\left[ 1 - \left(\frac{\mathcal{T}_{il}}{1 + \mathcal{T}_{il}}\right)^\beta \right] - \alpha_i\left( n_i - \sum_{l=1}^{m_i-1}(r_{il} +1) \right)\log\left[ 1 - \left(\frac{\mathcal{T}_{im_i}}{1 + \mathcal{T}_{im_i}}\right)^\beta \right].$$
 
 Now, it can be verified that the quantity 
 
 \begin{eqnarray} \nonumber
 \mathcal{P}_i(\beta) &=& -2\sum_{j=1}^{m_i-1}\log(\mathcal{Q}_{ij})\\ &=&\nonumber -2\sum_{j=1}^{m_i - 1}\log\left[  \frac{\sum_{l=1}^{j-1}(r_{il} + 1)\log\left[ 1 - \left(\frac{\mathcal{T}_{il}}{1 + \mathcal{T}_{il}}\right)^\beta \right] +                                                                                                                                                                                                                                                                                                                                                                                                                                                                                                                                                                                                                                                                                                                                                                                                                                                                                                                                                                                                                                                                                                                                                                                                                                                                                                                                                                                                                                                                                                                                                                                                                                                                                                                                                                                                                                                                                                                                                                                                                                                                                                                                                                                                                                                                                                                                                                                                                                                                                                                                                                                                                                                                                                                                                                                                                                                                                             \left( n_i - \sum_{l=1}^{j-1}(r_{il} +1) \right)\log\left[ 1 - \left(\frac{\mathcal{T}_{ij}}{1 + \mathcal{T}_{ij}}\right)^\beta \right]}{\sum_{l=1}^{m_i-1}(r_{il} + 1)\log\left[ 1 - \left(\frac{\mathcal{T}_{il}}{1 + \mathcal{T}_{il}}\right)^\beta \right] + \left( n_i - \sum_{l=1}^{m_i-1}(r_{il} +1) \right)\log\left[ 1 - \left(\frac{\mathcal{T}_{im_i}}{1 + \mathcal{T}_{im_i}}\right)^\beta \right]} \right]
 \end{eqnarray}
 has $\chi^2$ distribution with $2(m_i - 1)$ degrees of freedom, and so $\mathcal{P}_i(\beta)$ is a pivotal quantity since its distribution is independent of $\beta$. 
 
 Finally, the quantities $\mathcal{P}_i(\beta)$'s are mutually independent, and hence the quantity $\mathcal{P}(\beta) = \sum_{i=1}^{k}\mathcal{P}_i(\beta)$ has $\chi^2$ distribution with $2\sum_{i=1}^{k}(m_i - 1)$ degrees of freedom. The following theorem shows that the function $\mathcal{P}(\beta)$ is strictly decreasing when $\beta > 0$.

 \begin{theorem}\label{thm2}
 	The pivotal quantity $\mathcal{P}(\beta) = \sum_{i=1}^{k}\mathcal{P}_i(\beta)$ is a strictly decreasing function of $\beta > 0$.
 	
 	\begin{proof}
 	It should be observed that the reciprocal of the term written inside the logarithm in the expression of the pivotal quantity $\mathcal{P}_i(\beta)$ is as follows:
 	
 	\begin{eqnarray}\nonumber
 	 \left[  \frac{\sum_{l=1}^{m_i-1}(r_{il} + 1)\log\left[ 1 - \left(\frac{\mathcal{T}_{il}}{1 + \mathcal{T}_{il}}\right)^\beta \right] +                                                                                                                                                                                                                                                                                                                                                                                                                                                                                                                                                                                                                                                                                                                                                                                                                                                                                                                                                                                                                                                                                                                                                                                                                                                                                                                                                                                                                                                                                                                                                                                                                                                                                                                                                                                                                                                                                                                                                                                                                                                                                                                                                                                                                                                                                                                                                                                                                                                                                                                                                                                                                                                                                                                                                                                                                                                                                             \left( n_i - \sum_{l=1}^{m_i-1}(r_{il} +1) \right)\log\left[ 1 - \left(\frac{\mathcal{T}_{im_i}}{1 + \mathcal{T}_{im_i}}\right)^\beta \right]}{\sum_{l=1}^{j-1}(r_{il} + 1)\log\left[ 1 - \left(\frac{\mathcal{T}_{il}}{1 + \mathcal{T}_{il}}\right)^\beta \right] + \left( n_i - \sum_{l=1}^{j-1}(r_{il} +1) \right)\log\left[ 1 - \left(\frac{\mathcal{T}_{ij}}{1 + \mathcal{T}_{ij}}\right)^\beta \right]} \right] &=& \\ 1 +  \left[\frac{\sum_{l=j+1}^{m_i-1}(r_{il} + 1) \left(\frac{\phi_{il}(\beta)}{\phi_{ij}(\beta)}\right) + \left( n_i - \sum_{l=1}^{m_i-1}(r_{il} + 1) \right)\left(\frac{\phi_{im_i}(\beta)}{\phi_{ij}(\beta)}\right) - \left( n_i - \sum_{l=1}^{j}(r_{il} + 1) \right)}{ \sum_{l=1}^{j-1}(r_{il} + 1)\left(\frac{\phi_{il}(\beta)}{\phi_{ij}(\beta)}\right) + \left( n_i - \sum_{l=1}^{j-1}(r_{il} + 1) \right)} \right],
 	 \label{thm}
 	\end{eqnarray}
where $\phi_{ij}(\beta) = \log\left[ 1 - \left(\frac{\mathcal{T}_{ij}}{1 + \mathcal{T}_{ij}}\right)^\beta \right],~(j = 1,2, \dots, m_i - 1) .$  Now, we observe that $\frac{\phi_{il}(\beta)}{\phi_{ij}(\beta)}$ is an increasing function of $\beta$ if $l<j$ and a decreasing function of $\beta$ if $l > j$ because, $$\frac{d}{d\beta}\left[\frac{\phi_{il}(\beta)}{\phi_{ij}(\beta)}\right] = \frac{a_{ij}^{\beta}(1-a^{\beta}_{il})\log(a_{ij})\log(1 - a_{il}^{\beta}) - a_{il}^{\beta}(1 - a_{ij}^{\beta})\log(a_{il})\log(1 - a_{ij}^{\beta})}{(1 - a_{ij}^{\beta})(1-a^{\beta}_{il})\left[\log(1 - a_{ij}^{\beta})\right]^2},$$ where $a_{ij} = \frac{\mathcal{T}_{ij}}{1 + \mathcal{T}_{ij}},~(j = 1,2, \dots, m_i - 1)$, is positive if $l < j$ and negative if $l > j$. Therefore, in the right hand side of  (\ref{thm}), the numerator of the fraction inside the large parenthesis is a decreasing function of $\beta$ and the denominator of the same fraction is an increasing function of $\beta$. Also, it can be shown that, both the numerator and denominator are positive-valued functions of $\beta$. These facts finally imply that, the reciprocal term in (\ref{thm}) is a decreasing function of $\beta$. So, the pivotal quantity $\mathcal{P}_i(\beta)$ as well as the quantity $\mathcal{P}(\beta)$ is a strictly decreasing function of $\beta$.
    \end{proof}
 \end{theorem}
 Further, since $\mathcal{P}(\beta)$ is a continuous function of $\beta$, Theorem \ref{thm2} implies that for any $\beta^* > 0$, the equation $\mathcal{P}(\beta) = \beta^*$ has a unique solution $S(\beta^*)$, which can be considered as an estimate of $\beta$.\\
 On the other side, the following quantity
 \begin{equation}
 \mathcal{R}_i = 2\mathcal{W}_{im_i} = - 2\alpha_i\sum_{l=1}^{m_i-1}(r_{il} + 1)\log\left[ 1 - \left(\frac{\mathcal{T}_{il}}{1 + \mathcal{T}_{il}}\right)^\beta \right] - 2\alpha_i\left( n_i - \sum_{l=1}^{m_i-1}(r_{il} +1) \right)\log\left[ 1 - \left(\frac{\mathcal{T}_{im_i}}{1 + \mathcal{T}_{im_i}}\right)^\beta \right]
 \end{equation}
 has $\chi^2$ distribution with $2m_i$ degrees of freedom. When a pivotal estimate $S(\beta^*)$ of $\beta$ is obtained, an estimate $\alpha_i^*$ for $\alpha_i$ can be obtained as $ \alpha_i^* = \frac{\mathcal{R}^*_i}{2\Phi_i(S(\beta^*))}$, where the function $\Phi_i$ is defined as 
 \begin{equation}
 \Phi_i(\beta) =  - \sum_{l=1}^{m_i-1}(r_{il} + 1)\log\left[ 1 - \left(\frac{\mathcal{T}_{il}}{1 + \mathcal{T}_{il}}\right)^\beta \right] - \left( n_i - \sum_{l=1}^{m_i-1}(r_{il} +1) \right)\log\left[ 1 - \left(\frac{\mathcal{T}_{im_i}}{1 + \mathcal{T}_{im_i}}\right)^\beta \right]
 \end{equation}
and $\mathcal{R}_i^*$ is an unit sample drawn from the $\chi^2(2m_i)$ distribution. 
Furthermore, pivotal estimates for the unknown parameters $\alpha$ and $\beta$ as well as for the reliability and hazard functions, and MTF can be computed numerically through the Algorithm \ref{alg1} as follows:

\begin{algorithm}[H]
	\caption{Pivotal estimation for the RCs}
	\label{alg1}
	\begin{algorithmic}[1]
    \State Set $s = 1$.
    \State \label{s2} Generate a random sample {$\beta^*$} of unit size from $\chi^2$ distribution with $2\sum_{i=1}^{k}(m_i - 1)$ degrees of freedom. Let $\beta^{(s)}$ be the unique solution of the equation  $\mathcal{P}(\beta) = \beta^*$.
    \State Generate the random sample $\rho_i$ of unit size from the $\chi^2$ distribution with $2m_i$ degrees of freedom and then compute an estimate $\alpha_i^{(s)}$ for $\alpha_i$ as $\alpha_i^{(s)} = \frac{\rho_i}{2\Phi_i(\beta^{(s)})}$ for $i = 1, 2, \dots, k$.
    \State The pivotal estimate for $\alpha$ at the current step is computed as $\displaystyle \alpha^{(s)} = \left(\sum_{i=1}^{k} \frac{1}{\text{Var}(\alpha_i|\mathcal{T})} \alpha_i^{(s)}\right)/\left(\sum_{i=1}^{k} \frac{1}{\text{Var}(\alpha_i|\mathcal{T})}\right)$,  \text{where} $ \displaystyle \text{Var}(\alpha_i|\mathcal{T}) = \frac{1}{s}\sum_{l=1}^{s}\left( \alpha_i^{(l)} - \frac{1}{s} \sum_{l=1}^{s}\alpha_i^{(l)}\right)^2$.
    
    Further, the pivotal estimates for the reliability function, hazard function and the MTF for the current step are
    
    $$ \displaystyle{ R(t; \alpha^{(s)}, \beta^{(s)}) = \left[ 1 - \left( \frac{t}{1+t} \right)^{\beta^{(s)}} \right]^{\alpha^{(s)}}},$$
    
    $$ \displaystyle{ H(t; \alpha^{(s)}, \beta^{(s)}) = \frac{\alpha^{(s)}\beta^{(s)} t^{\beta^{(s)} - 1}(1 + t)^{-(\beta^{(s)} + 1)}}{1 - \left(\frac{t}{1+t} \right)^{\beta^{(s)}}}},$$
    
    and
    
    $$\displaystyle{ \mu(\alpha^{(s)}, \beta^{(s)}) =  \left\{(1 - 2^{-1/\alpha^{(s)}})^{-1/\beta^{(s)}} - 1 \right\}^{-1}  },$$
    
    respectively.
    
    \State $s = s + 1$ \label{s5}.
    
    \State Go through from the step \ref{s2} to step \ref{s5} repeatedly for a large number (say, $N$) of times to generate the large samples ($\theta^{(1)}, \theta^{(2)}, \dots, \theta^{(N)}$), where $\theta$ may represent any of the parameters $\alpha_1, \alpha_2, \dots, \alpha_k, \alpha, \beta$ or may represent the RCs $R(t; \alpha, \beta)$, $H(t; \alpha, \beta)$ or $\mu(\alpha, \beta)$. The pivotal estimate of $\theta$ is approximated as $\hat{\theta}_P = \frac{1}{N}\sum_{i=1}^{N}\theta^{(i)}$.
    
    \State The approximate variance of $\hat{\theta}_P$ is 
    
    $\mbox{Var}(\theta_P) = \frac{1}{N}\sum_{i=1}^{N}\left( \theta^{(i)} -  \hat{\theta}_P \right)^2$, where $\theta$ is the same thing as stated in the previous step.
    
    \State To compute the generalized confidence interval of $\theta$ (where $\theta$ may represent one of the parameters $\alpha_1, \alpha_2, \dots, \alpha_k, \beta$ or any one of the functions $R(t; \alpha, \beta)$, $H(t; \alpha, \beta)$ or $\mu(\alpha, \beta)$), arrange the samples $\theta^{(1)}, \theta^{(2)}, \dots, \theta^{(N)}$ in ascending order. If $\theta_{(1)}, \theta_{(2)}, \dots, \theta_{(N)}$ be the ordered sample, the approximate $100(1-\gamma)\%$ symmetric generalized comfidence interval for $\theta$ is $\left(\theta_{[N\gamma/2]}, \theta_{[N(1-\gamma/2)]}\right)$, where $[x]$ represents the largest integer less than or equal to $x$.
 	\end{algorithmic}
\end{algorithm}

\section{Order Statistics and Applications}\label{section5}
In this section, we look into a few distributional characteristics of the order statistics associated with the IEP distribution. Considering a random sample of size $n$, represented by $X_{1},X_{2},\ldots,X_{n}$, where $i=1,2,\ldots,n$. The random variables $X_{(1)},X_{(2)},\ldots,X_{(n)}$ are the order statistics. They are defined as function of $X_{1},X_{2},\ldots,X_{n}$, so that $X_{(1)}=\inf\{X_{1},X_{2},\ldots,X_{n}\}$, $X_{(n)}=\sup\{X_{1},X_{2},\ldots,X_{n}\}$, and $\mathbb{P}\{X_{(1)},X_{(2)},\ldots,X_{(n)}\}=1$. In theoretical research and practical applications, order statistics are important, particularly when it comes to the minimum, maximum, range, and analysis of a given $X_{i}$. The role of order statistics and records are found in \cite{barter1988history}. Order statistics are used in many different domains, including wireless communication, signal processing, quality control, survival analysis, life testing, reliability, and classification analysis. For more details on the concept of order statistics and its applications, please refer to \cite{david2004order}, \cite{govindaraju2013machine}, and \cite{shrahili2023excess}. Furthermore, order statistics are used in biomedical research (see \cite{greenberg1958applications}), image processing, and filtering theory. Additionally, order statistics filters, which are a type of non-linear filter (refer to \cite{dytso2021most}) and sampling plans (refer to \cite{schneider199818}).  

The CDF and PDF of the 1st or minimum order statistic of a random sample of size $n$ for the IEP distribution can be obtained as follows:
\begin{eqnarray}
F_{(1)}(t)=1-\left(1-\left(\frac{t}{1+t}\right)^{\beta}\right)^{n\alpha}~\mbox{and}~~f_{(1)}(t)=n\alpha\beta t^{\beta-1}(1+t)^{-(\beta+1)}\left(1-\left(\frac{t}{1+t}\right)^{\beta}\right)^{n\alpha-1}.
\end{eqnarray}
Therefore, the CDF and PDF of $r$th order statistic for the IEP distribution of $n$ size random sample are given as:
\begin{eqnarray}
F_{(r)}(t)=\sum_{j=r}^{n}\binom{n}{j}F^{j}(t)(1-F(t))^{n-j}=\sum_{j=r}^{n}\binom{n}{j}\left(1-\left(1-\left(\frac{t}{1+t}\right)^{\beta}\right)^{\alpha}\right)^{j} \left(1-\left(\frac{t}{1+t}\right)^{\beta}\right)^{\alpha(n-j)}  	
\end{eqnarray}
and 
\begin{eqnarray}
f_{(r)}(t)&=&\frac{n!}{(r-1)!\times1!\times(n-r)!}F^{r-1}(t)f(t)(1-F(t))^{n-r}\nonumber\\
&=&\frac{n!}{(r-1)!\times(n-r)!}\left(1-\left(1-\left(\frac{t}{1+t}\right)^{\beta}  \right)^{\alpha}\right)^{r-1}\alpha\beta t^{\beta-1}(1+t)^{-(\beta+1)}\left(1-\left(\frac{t}{1+t}\right)^{\beta}  \right)^{\alpha-1}\nonumber\\
&&\times \left(1-\left(\frac{t}{1+t}\right)^{\beta}\right)^{\alpha(n-r)}\nonumber\\
&=&\frac{n!\alpha\beta t^{\beta-1}(1+t)^{-(\beta+1)}}{(r-1)!\times(n-r)!}\left(1-\left(\frac{t}{1+t} \right)^{\beta}\right)^{\alpha(n-r+1)-1} \left(1-\left(1-\left(\frac{t}{1+t}\right)^{\beta}\right)^{\alpha}\right)^{r-1},   	
\end{eqnarray}
respectively. Finally, for the IEP distribution, the CDF and PDF of the $n$th or maximum order statistic of a random sample of size $n$ are as follows:
\begin{eqnarray}
F_{(n)}(x)=\left(1-\left(1-\left(\frac{t}{1+t}\right)^{\beta}\right)^{\alpha}  \right)^{n}  
\end{eqnarray}
and
\begin{eqnarray}
f_{(n)}(x)=n\alpha\beta t^{\beta-1}(1+t)^{-(\beta+1)} \left(1-\left(\frac{t}{1+t}\right)^{\beta}\right)^{\alpha-1}\left(1-\left(1-\left(\frac{t}{1+t}\right)^{\beta}\right)^{\alpha}\right)^{n-1}, 
\end{eqnarray}
respectively.

\section{Simulation Study} \label{sec 6}

This section is devoted for a detailed simulation study to assess the performances of the proposed estimation methods : the MLE method and the pivotal estimation method. At first, we manually constructed some BAPCS setups (listed within the Table \ref{t1}) with different numbers of total units and observable failures. Each BAPCS setup is utilized to simulate 2500 samples from the IEP($\alpha$, $\beta$) distribution with $\alpha = 3.5$ and $\beta = 2.25$. Now, using each of these samples, the MLE and pivotal estimates for $\alpha_i$ ($i=1,2,\cdots,k$), $\alpha$, and $\beta$ as well as for $R(t; \alpha, \beta)$, $H(t; \alpha, \beta)$, (taking $t = 0.75$) and $\mu(\alpha, \beta)$ and the generalized confidence intervals (GCIs) and the asymptotic confidence intervals (ACIs) for the same quantities are computed as shown in Sections \ref{sec 2}, \ref{sec 3}, and \ref{sec 4}. The estimated DDTFs can be visualized by observing the differences among the estimated values and estimated intervals of $\alpha_i$'s in different test facilities. To construct every one BAPCS setup, we implemented three different adaptive progressive Type-II censoring schemes manually and each of them is employed in every block of each BAPCS setup. Let $m$ be the total number of observable failures and $n$ be the total number of units under test in an arbitray block of an arbitray BAPCS setup. The three pre-determined censoring plans corresponding to the three adaptive progressive Type-II censoring schemes are then given as follows:

\begin{enumerate}
\item $R_1 = R_2 =  \dots = R_{\lfloor 3m/4 \rfloor-1} = 0, R_{\lfloor 3m/4 \rfloor } = \lceil \frac{n-m}{2} \rceil, R_{\lfloor 3m/4 \rfloor + 1} = \cdots = R_{m-1} = 0, R_{m} = n  - m - \lceil \frac{n-m}{2} \rceil $.
\item $R_1 = R_2 =  \dots = R_{\lfloor m/4 \rfloor-1} = 0, R_{\lfloor m/4 \rfloor } = \lceil \frac{n-m}{2} \rceil, R_{\lfloor m/4 \rfloor + 1 }  = \cdots = R_{\lfloor 3m/4 \rfloor - 1} = 0, R_{\lfloor 3m/4 \rfloor} = n  - m - \lceil \frac{n-m}{2} \rceil, R_{\lfloor 3m/4 \rfloor + 1} = \cdots = R_m = 0$.
\item $R_1 = R_2 =  \dots = R_{\lfloor m/4 \rfloor-1} = 0, R_{\lfloor m/4 \rfloor } = \lceil \frac{n-m}{2} \rceil, R_{\lfloor m/4 \rfloor + 1 }  = \cdots = R_{m-1} = 0, R_m = n  - m - \lceil \frac{n-m}{2} \rceil$.
\end{enumerate}

Here, the notations $\lfloor . \rfloor$ and $\lceil . \rceil$ represent the floor function and the ceiling function, respectively. For each BAPCS setup, to maintain simplicity, the threshold times $T_i$'s are taken equal for every block or facility.
 
All the estimated values, approximate biases and variances , ACIs, GCIs, and the lengths (denoted by L in the Tables \ref{t2} to \ref{t7}) of ACIs and GCIs of $\alpha$, $\beta$ as well as of $R(t; \alpha, \beta)$, $H(t; \alpha, \beta)$ (taking $t = 0.75$), and $\mu( \alpha, \beta)$ corressponding to each BAPCS setup and each pe-assumed censoring plan (mentioned under the column CP of the Tables \ref{t2}  to \ref{t7}) are listed from Table \ref{t2} to Table \ref{t7}. 
From Table \ref{t2} to Table \ref{t7}, the following observations are readily available:

\begin{itemize}
\item[i.] In most of the cases, absolute bias and variance of each estimator are smaller in case of larger values of the number of total testing units $n$.
\item[ii.] Lengths of ACIs and GCIs become smaller in case of larger values of $n$. 
\item[iii.] In most of the cases, pivotal estimation produces smaller absolute biases and variances compared to the MLE method.
\item[iv.] The GCIs are narrower than the ACIs in majority of the cases.
\item[v.] In most of the situations, the absolute biases, variances, and ACI and GCI lengths of the estimators are smaller in case of the censoring plans 2 and 3 compared to the plan 1.
\end{itemize}

\begin{table}[H]
	\caption{\label{t1}\normalsize BAPCS setups for simulation study.}
	\setlength\tabcolsep{1pt}
	\begin{tabular}{|c|c|c|c|c|c|}
		\hline
		\thead{Setup\\No.}&\thead{Number of \\ total\\ units ($n$)} &\thead{Number of \\ blocks ($k$)} & \thead{$T_i$~ \\($i = 1,2,\cdots,k)$} & \thead{ Sample sizes \\in different \\ blocks \\($n_1, n_2, \dots, n_k$)}  & \thead{ Observable failures \\in different blocks \\ ($m_1,m_2, \dots, m_k$)}  \\
		\hline
		1&\thead{$200$}& \thead{$4$}  &$0.75$& \thead{$ S_1 = \{55,~45,~46,~54\} $}  & \thead{$E_1 = \{45,~36,~34,~45\}$ } \\
		\hline
		2&\thead{$225$}& \thead{$4$} & $0.75$ & \thead{$S_2 = \{55,~ 60,~50,~60\}$}  & \thead{ $E_2 = \{ 44,~50,~40,~46\}$ }\\
		\hline
		3&\thead{$250$}&\thead{$4$}&$0.75$&\thead{$S_3 = \{60,~65,~,65,~60\}$}&\thead{$E_3 = \{53,~57,~58,~52\}$}\\
		\hline
		4&\thead{$200$}& \thead{$5$}  &$0.5$& \thead{$ S_4 = \{38,~42,~43,~37,~40\} $}  & \thead{$E_4 = \{32,~32,~38,~28,~30\}$ } \\
		\hline
		5&\thead{$225$}& \thead{$5$}&0.5  & \thead{$ S_5 = \{44,~48,~40,~45,~48\} $}  & \thead{$E_5 = \{33,~39,~32,~36,~40\}$ } \\
		\hline
		6&\thead{$250$}& \thead{$5$}&0.5 & \thead{$ S_6 = \{50,~57,~45,~50, ~48\} $}  & \thead{$E_6 = \{45,~51, ~39,~45,~40\}$ } \\
		\hline
		
	\end{tabular}
\end{table}

\begin{table}[H]
	\centering
	\small
	\caption{\small Point and interval estimates under the BAPCS setup 1.}
	\label{t2}
	\setlength\tabcolsep{1.5pt}
	\renewcommand{\arraystretch}{1.2}
	\begin{tabular}{cccccccccccccccccc}
		\hline
		&&&  \multicolumn{6}{c}{ML Estimates} & \multicolumn{6}{c}{~~ Pivotal Estimates} \\
		\cmidrule(lr){4-9}  \cmidrule(lr){10-17}
		&&&Est.&Bias& Variance &\multicolumn{2}{c}{ACIs}&L&    Est.&Bias& Variance &\multicolumn{2}{c}{~~ ~~GCIs~~}&L&&   \\
		\cmidrule(lr){7-8}  \cmidrule(lr){13-14}
		$n$&CP&Pars.&&&  &Lower & Upper  & & & &  & Lower &Upper\\ \hline
		200&1& $ \beta$& 2.3120& 0.0620& 0.0353& 1.9445& 2.6795& 0.7349& 2.2645& 0.0145& 0.0350& 1.9073& 2.6380& 0.7307   \\
		&&$ \alpha_{1}$&
		3.7673& 0.2673& 0.6031& 2.2920& 5.2427& 2.9507& 3.6668& 0.1668& 0.5823& 2.4079&5.2916&2.8837  \\
		&& $ \alpha_{2}$&3.8000& 0.3000& 0.7166& 2.1980& 5.4019& 3.2039& 3.6973& 0.1973& 0.6915& 2.3426& 5.4680& 3.1254  \\
		&&$ \alpha_{3}$&3.8222& 0.3222& 0.7818& 2.1510& 5.4935& 3.3425& 3.7168& 0.2168& 0.7553& 2.3130& 5.5751& 3.2621  \\
		&&$ \alpha_{4}$&3.7666& 0.2666& 0.5944& 2.2990& 5.2342& 2.9352& 3.6667& 0.1667& 0.5730& 2.4142& 5.2779& 2.8637  \\
		&&$ \alpha$& 3.6143& 0.1143& 0.3257& 2.5177& 4.7108& 2.1931& 3.5181& 0.0181& 0.3126& 2.5623& 4.7010& 2.1387  \\
		&& $R(t)$&0.5801& 0.0107& 0.0008& 0.5244& 0.6359& 0.1114& 0.5776& 0.0082& 0.0008& 0.5207& 0.6327& 0.1120  \\
		&& $h(t)$ &1.0355& -0.0118& 0.0079& 0.8625& 1.2085& 0.3459& 1.0256& -0.0217& 0.0078& 0.8596& 1.2043& 0.3447  \\
		& &$\text{MTF}$ &0.8952& 0.0215& 0.0032& 0.7851& 1.0052& 0.2201& 0.8952& 0.0215& 0.0034& 0.7869& 1.0143& 0.2274 \\
		200&2& $\beta$& 2.3062& 0.0562& 0.0325& 1.9537& 2.6588& 0.7051& 2.2628& 0.0128& 0.0323& 1.9188& 2.6205& 0.7017\\
		&&$ \alpha_{1}$&3.7326& 0.2326& 0.5495& 2.3196& 5.1455&2.8259& 3.6432& 0.1432& 0.5325& 2.4235& 5.1915& 2.7680   \\
		&&$ \alpha_{2}$&3.7713& 0.2713& 0.6563& 2.2328& 5.4490&3.0770& 3.6804& 0.1804& 0.6362& 2.3633& 5.3729& 3.0096 \\
		&&$ \alpha_{3}$&3.7815& 0.2815& 0.7061& 2.1858&5.3772&3.1914&3.6880&0.1880&0.6831&2.3310&5.4490 &3.1180\\
		&&$ \alpha_{4}$&3.7356& 0.2356& 0.5445& 2.3285&5.1428&2.8143&3.6472&0.1474&0.5275&2.4327&5.1865 &2.7539\\
		&&$\alpha$&3.5872& 0.0872& 0.2828& 2.5629& 4.6115& 2.0486& 3.5019& 0.0019&0.2722& 2.5962& 4.5999& 2.0037 \\
		&& $R(t)$&0.5805& 0.0110& 0.0008& 0.5238& 0.6372&  0.1134& 0.5778& 0.0084& 0.0008& 0.5200&0.6399& 0.1139  \\
		&& $h(t)$ & 1.0323& -0.0150& 0.0075& 0.8637& 1.2009&0.3372& 1.0242& -0.0231& 0.0075& 0.8620& 1.1985& 0.3365  \\
		& &$\text{MTF}$ &0.8963& 0.0226& 0.0033& 0.7842& 1.0083& 0.2241& 0.8958& 0.0221& 0.0035& 0.7859&1.0164&0.2305  \\
		200&3& $ \beta$&2.3067& 0.0567& 0.0325& 1.9540& 2.6594&0.7054& 2.2631& 0.0131& 0.0323& 1.9188& 2.6216& 0.7028   \\
		&&$ \alpha_{1}$&3.7442& 0.2442& 0.5539& 2.3256& 5.1629& 2.8373&3.6542& 0.1542& 0.5364& 2.4295& 5.2092& 2.7797   \\
		&& $ \alpha_{2}$&3.7694& 0.2694& 0.6571& 2.2319& 5.3069&3.0750& 3.6781& 0.1781& 0.6369& 2.3617& 5.3727& 3.0110 \\
		& &$ \alpha_{3}$&3.7946& 0.2946& 0.7150& 2.1914& 5.3978& 3.2063& 3.7001& 0.2001& 0.6934& 2.3361& 5.4740& 3.1378  \\
		& &$ \alpha_{4}$&3.7550& 0.2550& 0.5512& 2.3389& 5.1711& 2.8322& 3.6654& 0.1654& 0.5336& 2.4433& 5.2154& 2.7721  \\
		&&$\alpha$& 3.5982& 0.0982& 0.2856& 2.5693& 4.6270& 2.0578& 3.5124& 0.0124& 0.2755& 2.6017& 4.6165& 2.0148   \\
		&& $R(t)$&0.5797& 0.0103& 0.0008& 0.5230& 0.6364& 0.1134& 0.5770& 0.0076& 0.0008& 0.5192& 0.6331& 0.1140   \\
		&& $h(t)$ &1.0350& -0.0123& 0.0075& 0.8658& 1.2041& 0.3383& 1.0270& -0.0204& 0.0075& 0.8643& 1.2016& 0.3373   \\
		& &$\text{MTF}$& 0.8947& 0.0211& 0.0033& 0.7830& 1.0065& 0.2236& 0.8941& 0.0205& 0.0035& 0.7845& 1.0145& 0.2300  \\
		\hline
	\end{tabular}
\end{table}

\begin{table}[H]
	\centering
	\small
	\caption{\small Point and interval estimates under the BAPCS setup 2.}
	\label{t3}
	\setlength\tabcolsep{1.5pt}
	\renewcommand{\arraystretch}{1.2}
	\begin{tabular}{cccccccccccccccccc}
		\hline
		&&&  \multicolumn{6}{c}{ML Estimates} & \multicolumn{6}{c}{~~ Pivotal Estimates} \\
		\cmidrule(lr){4-9}  \cmidrule(lr){10-17}
		&&&Est.&Bias& Variance &\multicolumn{2}{c}{ACIs}&L&    Est.&Bias& Variance &\multicolumn{2}{c}{~~ ~~GCIs~~}&L&&   \\
		\cmidrule(lr){7-8}  \cmidrule(lr){13-14}
		$n$&CP&Pars.&&&  &Lower & Upper  & & & &  & Lower &Upper\\ \hline
		225&1& $ \beta$& 2.3068& 0.0568& 0.0313& 1.9607& 2.6530& 0.6923& 2.2643& 0.0143& 0.0310& 1.9274& 2.6155& 0.6881  \\
		&&$ \alpha_{1}$&3.7559& 0.2559& 0.5809& 2.3035& 5.2083& 2.9048& 3.6648& 0.1648& 0.5617& 2.4208&5.2595&2.8388  \\
		&& $ \alpha_{2}$&3.7464& 0.2464& 0.5205& 2.3669& 5.1259& 2.7590& 3.6575& 0.1575& 0.5034& 2.4667& 5.1661& 2.6995  \\
		&&$ \alpha_{3}$&3.7560& 0.2560& 0.6190& 2.2608& 5.2512& 2.9904& 3.6655& 0.1655& 0.5980& 2.3872& 5.3094& 2.9223  \\
		&&$ \alpha_{4}$&3.7448& 0.2448& 0.5774& 2.2998& 5.1899& 2.8901& 3.6533& 0.1533& 0.5591& 2.4159& 5.2450& 2.8292  \\
		&&$ \alpha$& 3.5993& 0.0993& 0.2847& 2.5714& 4.6272& 2.0558& 3.5135& 0.0135& 0.2736& 2.6088& 4.6167& 2.0079  \\
		&& $R(t)$&0.5796& 0.0102& 0.0007& 0.5271& 0.6321& 0.1050& 0.5773& 0.0079& 0.0007& 0.5238& 0.6294& 0.1055  \\
		&& $h(t)$ &1.0352& -0.0121& 0.0070& 0.8725& 1.1980& 0.3255& 1.0264& -0.0209& 0.0069& 0.8699& 1.1939& 0.3240 \\
		& &$\text{MTF}$ &0.8939& 0.0202& 0.0028& 0.7904& 0.9974& 0.2071& 0.8938& 0.0201& 0.0030& 0.7921& 1.0053& 0.2132 \\
		225&2& $ \beta$& 2.3029& 0.0529& 0.0289& 1.9705& 2.6354& 0.6648& 2.2639& 0.0139& 0.0286& 1.9394& 2.6005& 0.6611  \\
		&&$ \alpha_{1}$&
		3.7182& 0.2182& 0.5296& 2.3282& 5.1083& 2.7801& 3.6371& 0.1371& 0.5138& 2.4311&5.1563&2.7252  \\
		&& $ \alpha_{2}$&3.7284& 0.2284& 0.4856& 2.3973& 5.0595& 2.6622& 3.6479& 0.1479& 0.4712& 2.4875& 5.0978& 2.6103  \\
		&&$ \alpha_{3}$&3.7763& 0.2763& 0.5890& 2.3148& 5.2377& 2.9229& 3.6933& 0.1963& 0.5714& 2.4315& 5.2940& 2.8625  \\
		&&$ \alpha_{4}$&3.7425& 0.2425& 0.5332& 2.3486& 5.1363& 2.7877& 3.6598& 0.1598& 0.5178& 2.4523& 5.1840& 2.7316  \\
		&&$ \alpha$& 3.5899& 0.0899& 0.2514& 2.6222& 4.5577& 1.9355& 3.5129& 0.0129& 0.2426& 2.6497& 4.5460& 1.8963  \\
		&& $R(t)$&0.5790& 0.0095& 0.0007& 0.5256& 0.6324& 0.1068& 0.5766& 0.0072& 0.0008& 0.5222& 0.6294& 0.1072  \\
		&& $h(t)$ &1.0359& -0.0114& 0.0067& 0.8764& 1.1954& 0.3189& 1.0286& -0.0187& 0.0067& 0.8747& 1.1929& 0.3182 \\
		& &$\text{MTF}$ &0.8929& 0.0192& 0.0029& 0.7879& 0.9978& 0.2099& 0.8923& 0.0186& 0.0030& 0.7893& 1.0046& 0.2153 \\
		225&3& $\beta$& 2.3016& 0.0516& 0.0289& 1.9692& 2.6340& 0.6648& 2.2624& 0.0124& 0.0286& 1.9380& 2.5989& 0.6608\\
		&&$ \alpha_{1}$&3.7158& 0.2158& 0.5297& 2.3267& 5.1050&2.7782& 3.6342& 0.1342& 0.5137& 2.4299& 5.1511& 2.7212   \\
		&&$ \alpha_{2}$&3.7272& 0.2272& 0.4861& 2.3962& 5.0581&2.6619& 3.6465& 0.1465& 0.4919& 2.4861& 5.0968& 2.6107 \\
		&&$ \alpha_{3}$&3.7567& 0.2567& 0.5818& 2.3036&5.2098&2.9062&3.6738&0.1738&0.5638&2.4192&5.2639 &2.8447\\
		&&$ \alpha_{4}$&3.7297& 0.2297& 0.5287& 2.3411&5.1184&2.7773&3.6471&0.1471&0.5135&2.4442&5.1667 &2.7224\\
		&&$\alpha$&3.5825& 0.0825& 0.2504& 2.6170& 4.5480& 1.9310& 3.5050& 0.0050&0.2415& 2.6444& 4.5348& 1.8904 \\
		&& $R(t)$&0.5793& 0.0098& 0.0007& 0.5259& 0.6326&  0.1068& 0.5769& 0.0075& 0.0007& 0.5225&0.6297& 0.1072  \\
		&& $h(t)$ & 1.0344& 0.0129& 0.0067& 0.8752& 1.1937&0.3184& 1.0271& -0.0202& 0.0066& 0.8734& 1.1911& 0.3177  \\
		& &$\text{MTF}$ &0.8935& 0.0199& 0.0029& 0.7884& 0.9987& 0.2102& 0.8930& 0.0194& 0.0031& 0.7899&1.0055&0.2157  \\
		\hline
	\end{tabular}
\end{table}

\begin{table}[H]
	\centering
	\small
	\caption{\small Point and interval estimates under the BAPCS setup 3.}
	\label{t4}
	\setlength\tabcolsep{1.5pt}
	\renewcommand{\arraystretch}{1.2}
	\begin{tabular}{cccccccccccccccccc}
		\hline
		&&&  \multicolumn{6}{c}{ML Estimates} & \multicolumn{6}{c}{~~ Pivotal Estimates} \\
		\cmidrule(lr){4-9}  \cmidrule(lr){10-17}
		&&&Est.&Bias& Variance &\multicolumn{2}{c}{ACIs}&L&    Est.&Bias& Variance &\multicolumn{2}{c}{~~ ~~GCIs~~}&L&&   \\
		\cmidrule(lr){7-8}  \cmidrule(lr){13-14}
		$n$&CP&Pars.&&&  &Lower & Upper  & & & &  & Lower &Upper\\ \hline
		250&1& $ \beta$& 2.2939& 0.0439& 0.0251& 1.9840& 2.6038& 0.6167& 2.2604& 0.0104& 0.0249& 1.9575& 2.5741& 0.6167   \\
		&&$ \alpha_{1}$&3.6826& 0.1826& 0.4267& 2.4291& 4.9362& 2.5071& 3.6156& 0.1156& 0.4157& 2.5131&4.9759&2.4628  \\
		&& $ \alpha_{2}$&3.7047& 0.2047& 0.4159& 2.4682& 4.9412& 2.4731& 3.6369& 0.1369& 0.4054& 2.5475& 4.9788& 2.4313  \\
		&&$ \alpha_{3}$&3.6746& 0.1746& 0.3994& 2.4628& 4.8865& 2.4237& 3.6081& 0.1081& 0.3902& 2.5384& 4.9241& 2.3830  \\
		&&$ \alpha_{4}$&3.6988& 0.1988& 0.4427& 2.4248& 4.9728& 2.5479& 3.6303& 0.1303& 0.4310& 2.5132& 5.0157& 2.5025  \\
		&&$ \alpha$& 3.5660& 0.0660& 0.2037& 2.6926& 4.4394& 1.7468& 3.5011& 0.0011& 0.1975& 2.7152& 4.4289& 1.7136  \\
		&& $R(t)$&0.5781& 0.0087& 0.0006& 0.5288& 0.6275& 0.0987& 0.5759& 0.0065& 0.0006& 0.5256& 0.6248& 0.0993  \\
		&& $h(t)$ &1.0351& -0.0122& 0.0053& 0.8922& 1.1780& 0.2858& 1.0292& -0.0181& 0.0053& 0.8908& 1.1761& 0.2853  \\
		& &$\text{MTF}$ &0.8910& 0.0173& 0.0024& 0.7948& 0.9872& 0.1924& 0.8900& 0.0163& 0.0025& 0.7950& 0.9918& 0.1968 \\
		250&2& $\beta$& 2.2914& 0.0414& 0.0237& 1.9899& 2.5929& 0.6030& 2.2600& 0.0100& 0.0236& 1.9642& 2.5651& 0.6008\\
		&&$ \alpha_{1}$&3.6969& 0.1969& 0.4134& 2.4629& 4.9309&2.4680& 3.6345& 0.1345& 0.4044& 2.5425&4.9740& 2.4315   \\
		&&$ \alpha_{2}$&3.6740& 0.1740& 0.3910& 2.4742& 4.8739&2.3997& 3.6112& 0.1112& 0.3821& 2.5478& 4.9123& 2.3645 \\
		&&$ \alpha_{3}$&3.6616& 0.1616& 0.3787& 2.4792&4.8441&2.3649&3.5999&0.0999&0.3700&2.5512&4.8785 &2.3273\\
		&&$ \alpha_{4}$&3.6831& 0.1831& 0.4200& 2.4402&4.9261&2.4859&3.6202&0.1202&0.4100&2.5232&4.9682 &2.4449\\
		&&$\alpha$&3.5598& 0.0598& 0.1876& 2.7209& 4.3988& 1.6778& 3.5000& 0.0000&0.1824& 2.7384& 4.3889& 1.6504 \\
		&& $R(t)$&0.5778& 0.0083& 0.0006& 0.5279& 0.6277&  0.0998& 0.5756& 0.0061& 0.0007& 0.5248&0.6250& 0.1002  \\
		&& $h(t)$ & 1.0354& -0.0119& 0.0052& 0.8940& 1.1769&0.2829& 1.0304& -0.0169& 0.0052& 0.8832& 1.1757& 0.2825  \\
		& &$\text{MTF}$ &0.8906& 0.0170& 0.0025& 0.7934& 0.9878& 0.1944& 0.8895& 0.0158& 0.0026& 0.7938&0.9921&0.1983  \\
		250&3& $ \beta$&2.2881& 0.0381& 0.0237& 1.9869& 2.5893&0.6024& 2.2569& 0.0069& 0.0236& 1.9614& 2.5615& 0.6001   \\
		&&$ \alpha_{1}$&3.6812& 0.1812& 0.4090& 2.4531& 4.9093& 2.4562&3.6186& 0.1186& 0.3996& 2.5320& 4.9495& 2.4174   \\
		&& $ \alpha_{2}$&3.6660& 0.1660& 0.3871& 2.4696& 4.8625&2.3929& 3.6047& 0.1047& 0.3787& 2.5435& 4.9002& 2.3567 \\
		& &$ \alpha_{3}$&3.6515& 0.1515& 0.3758& 2.4726& 4.8303& 2.3577& 3.5900& 0.0900& 0.3678& 2.5426& 4.8648& 2.3222  \\
		& &$ \alpha_{4}$&3.6736& 0.1736& 0.4162& 2.4345& 4.9127& 2.4782& 3.6112& 0.1112& 0.4066& 2.5171& 4.9558& 2.4387  \\
		&&$\alpha$& 3.5514& 0.0514& 0.1864& 2.7149& 4.3878& 1.6729& 3.4920& -0.008& 0.1814& 2.7312& 4.3771& 1.6459   \\
		&& $R(t)$&0.0.5775& 0.0081& 0.0006& 0.5277& 0.6274& 0.0998& 0.5753& 0.0059& 0.0007& 0.5246& 0.6248& 0.1002   \\
		&& $h(t)$ &1.0350& -0.0123& 0.0052& 0.8936& 1.1763& 0.2827& 1.0299& -0.0174& 0.0052& 0.8926& 1.1753& 0.2826   \\
		& &$\text{MTF}$& 0.8902& 0.0166& 0.0025& 0.7930& 0.9875& 0.1945& 0.8891& 0.0154& 0.0026& 0.7934& 0.9918& 0.1984  \\
		\hline
	\end{tabular}
\end{table}

\begin{table}[H]
	\centering
	\small
	\caption{\small  Point and interval estimates under the BAPCS setup 4.}
	\label{t5}
	\setlength\tabcolsep{1.5pt}
	\renewcommand{\arraystretch}{1.2}
	\begin{tabular}{cccccccccccccccccc}
		\hline
		&&&  \multicolumn{6}{c}{ML Estimates} & \multicolumn{6}{c}{~~ Pivotal Estimates} \\
		\cmidrule(lr){4-9}  \cmidrule(lr){10-17}
		&&&Est.&Bias& Variance &\multicolumn{2}{c}{ACIs}&L&    Est.&Bias& Variance &\multicolumn{2}{c}{~~ ~~GCIs~~}&L&&   \\
		\cmidrule(lr){7-8}  \cmidrule(lr){13-14}
		$n$&CP&Pars.&&&  &Lower & Upper  & & & &  & Lower &Upper\\ \hline
		200&1& $ \beta$& 2.3199& 0.0699& 0.0356& 1.9511& 2.6887& 0.7376& 2.2601& 0.0101& 0.0351& 1.9022& 2.6341& 0.7319&   \\
		&&$ \alpha_{1}$&3.8099& 0.3099& 0.7606& 2.1629& 5.4569& 3.294& 3.6785& 0.1785& 0.7214& 2.2985& 5.4846& 3.1861&  \\
		&& $ \alpha_{2}$&3.8373& 0.3373& 0.8107& 2.1376& 5.5369& 3.3993& 3.699& 0.199& 0.7679& 2.2856& 5.5676& 3.282& 
		\\
		&&$ \alpha_{3}$& 3.7801& 0.2801& 0.6448& 2.2540& 5.3062& 3.0523& 3.6518& 0.1518& 0.611& 2.3614& 5.3158& 2.9543& \\
		&&$ \alpha_{4}$&3.8706& 0.3706& 0.9018& 2.0770& 5.6643& 3.5873& 3.7298& 0.2298& 0.8538& 2.2463& 5.7129& 3.4666&  \\
		&&$ \alpha_{5}$&3.8589& 0.3589& 0.8749& 2.1053& 5.6124& 3.5071& 3.7183& 0.2183& 0.8285& 2.2657& 5.6514& 3.3858&  \\
		&&$ \alpha$& 3.5988& 0.0988& 0.3221& 2.5071& 4.6906& 2.1835& 3.4745& -0.0255& 0.3033& 2.5315& 4.6397& 2.1082&  \\
		&& $R(t)$&0.5836& 0.0142& 0.0008& 0.5278& 0.6394& 0.1116& 0.5801& 0.0107& 0.0008& 0.523& 0.6354& 0.1124&  \\
		&& $h(t)$ & 1.0272& -0.0201& 0.0078& 0.8554& 1.1991& 0.3437& 1.0160& -0.0313& 0.0077& 0.8514& 1.1932& 0.3418&  \\
		& &$\text{MTF}$ &0.9023& 0.0286& 0.0033& 0.7910& 1.0136& 0.2227& 0.9011& 0.0274& 0.0035& 0.7913& 1.022& 0.2307& \\
		200&2& $ \beta$& 2.3139& 0.0639& 0.0330& 1.9588& 2.669& 0.7102& 2.2615& 0.0115& 0.0327& 1.9156& 2.6218& 0.7061   \\
		&&$ \alpha_{1}$&3.7918& 0.2918& 0.7122& 2.1895& 5.3942& 3.2047& 3.6784& 0.1784& 0.6822& 2.3182& 5.4363& 3.118&  \\
		&& $ \alpha_{2}$&3.8263& 0.3263& 0.7566& 2.1800& 5.4726& 3.2926& 3.7087& 0.2087& 0.7257& 2.3195& 5.5234& 3.2039& 
		\\
		&&$ \alpha_{3}$& 3.7659& 0.2659& 0.6070& 2.2813& 5.2506& 2.9693& 3.6564& 0.1564& 0.5823& 2.3849& 5.2764& 2.8915 \\
		&&$ \alpha_{4}$&3.8062& 0.3062& 0.8217& 2.0951& 5.5174& 3.4223& 3.6895& 0.1895& 0.7846& 2.2526& 5.575& 3.3224&  \\
		&&$ \alpha_{5}$&3.8256& 0.3256& 0.8001& 2.1398& 5.5114& 3.3716& 3.7068& 0.2068& 0.7645& 2.2904& 5.5644& 3.2739&  \\
		&&$ \alpha$& 3.5835& 0.0835& 0.2863& 2.5533& 4.6137& 2.0605& 3.4771& -0.0229& 0.2727& 2.5726& 4.5762& 2.0036&  \\
		&& $R(t)$&0.5831& 0.0137& 0.0008& 0.5265& 0.6398& 0.1133& 0.5798& 0.0104& 0.0008& 0.5219& 0.6359& 0.1141&  \\
		&& $h(t)$ & 1.0265& -0.0208& 0.0074& 0.8583& 1.1948& 0.3365& 1.0174& -0.0300& 0.0074& 0.8557& 1.1913& 0.3356&  \\
		& &$\text{MTF}$ &0.9014& 0.0278& 0.0033& 0.7887& 1.0142& 0.2256& 0.9001& 0.0265& 0.0036& 0.7891& 1.0219& 0.2327& \\
		200&3& $ \beta$& 2.3118& 0.0618& 0.0329& 1.9569& 2.6666& 0.7097& 2.2593& 0.0093& 0.0326& 1.9137& 2.6192& 0.7055   \\
		&&$ \alpha_{1}$&3.7768& 0.2768& 0.7081& 2.1809& 5.3726& 3.1917& 3.6633& 0.1633& 0.6766& 2.3117& 5.4122& 3.1005&  \\
		&& $ \alpha_{2}$&3.8298& 0.3298& 0.7597& 2.1807& 5.4789& 3.2982& 3.7112& 0.2112& 0.7263& 2.3209& 5.5211& 3.2002& 
		\\
		&&$ \alpha_{3}$& 3.7813& 0.2813& 0.6135& 2.2895& 5.2731& 2.9836& 3.6716& 0.1716& 0.5883& 2.395& 5.3001& 2.9051& \\
		&&$ \alpha_{4}$&3.8287& 0.3287& 0.8327& 2.1061& 5.5513& 3.4452& 3.7101& 0.2101& 0.7960& 2.2654& 5.6098& 3.3444&  \\
		&&$ \alpha_{5}$&3.8187& 0.3187& 0.7913& 2.1374& 5.5000& 3.3626& 3.6998& 0.1998& 0.7570& 2.2866& 5.5540& 3.2674&  \\
		&&$ \alpha$& 3.5870& 0.0870& 0.2873& 2.5551& 4.6189& 2.0639& 3.4805& -0.0195& 0.2736& 2.5736& 4.5816& 2.0080&  \\
		&& $R(t)$&0.5821& 0.0127& 0.0008& 0.5254& 0.6388& 0.1134& 0.5788& 0.0094& 0.0008& 0.5208& 0.6349& 0.1141&  \\
		&& $h(t)$ & 1.0289& -0.0184& 0.0075& 0.8601& 1.1977& 0.3376& 1.0198& -0.0276& 0.0075& 0.8574& 1.1942& 0.3368&  \\
		& &$\text{MTF}$ &0.8995& 0.0258& 0.0033& 0.7869& 1.012& 0.225& 0.8981& 0.0244& 0.0035& 0.7874& 1.0197& 0.2323& \\
		\hline
	\end{tabular}
\end{table}

\begin{table}[H]
	\centering
	\small
	\caption{\small Point and interval estimates under the BAPCS setup 5.}
	\label{t6}
	\setlength\tabcolsep{1.5pt}
	\renewcommand{\arraystretch}{1.2}
	\begin{tabular}{cccccccccccccccccc}
		\hline
		&&&  \multicolumn{6}{c}{ML Estimates} & \multicolumn{6}{c}{~~ Pivotal Estimates} \\
		\cmidrule(lr){4-9}  \cmidrule(lr){10-17}
		&&&Est.&Bias& Variance &\multicolumn{2}{c}{ACIs}&L&    Est.&Bias& Variance &\multicolumn{2}{c}{~~ ~~GCIs~~}&L&&   \\
		\cmidrule(lr){7-8}  \cmidrule(lr){13-14}
		$n$&CP&Pars.&&&  &Lower & Upper  & & & &  & Lower &Upper\\ \hline
		225&1& $ \beta$& 2.3203& 0.0703& 0.0316& 1.9725& 2.6682& 0.6957& 2.2671& 0.0171& 0.0313& 1.9285& 2.6194& 0.6909   \\
		&&$ \alpha_{1}$&
		3.8305& 0.3305& 0.7592& 2.1818& 5.4793& 3.2975& 3.7063& 0.2063& 0.7227& 2.3221& 5.5204& 3.1983  \\
		&& $ \alpha_{2}$&
		3.8117& 0.3117& 0.6479& 2.2845& 5.3388& 3.0543& 3.6921& 0.1921& 0.6180& 2.3949& 5.3622& 2.9673  \\
		&&$ \alpha_{3}$&3.8188& 0.3188& 0.7501& 2.1819& 5.4558& 3.2740& 3.6972& 0.1972& 0.7140& 2.3193& 5.4932& 3.1739  \\
		&&$ \alpha_{4}$&3.8156& 0.3156& 0.6869& 2.4323& 5.3889& 3.1466& 3.6954& 0.1954& 0.6545& 2.3654& 5.4177& 3.0523  \\
		&&$\alpha_{5}$&3.8225& 0.3225& 0.6314& 2.3123& 5.3326& 3.0203& 3.7038& 0.2038& 0.6028& 2.4179& 5.3523& 2.9343\\
		&&$ \alpha$& 3.6189& 0.1189& 0.2906& 2.5817& 4.6562& 2.0745& 3.5069& 0.0069& 0.2753& 2.6007& 4.6135& 2.0127  \\
		&& $R(t)$&0.5820& 0.0125& 0.0007& 0.5294& 0.6345& 0.1051& 0.5790& 0.0095& 0.0007& 0.5252& 0.6310& 0.1058 \\
		&& $h(t)$ &1.0327& -0.0146& 0.0070& 0.8700& 1.1954& 0.3254& 1.0225& -0.0248& 0.0069& 0.8659& 1.1903& 0.3244  \\
		& &$\text{MTF}$ &0.8983& 0.0247& 0.0028& 0.7943& 1.0024& 0.2081& 0.8973& 0.0236& 0.0030& 0.7946& 1.0095& 0.2148 \\
		225&2& $ \beta$& 2.3073& 0.0573& 0.0291& 1.9736& 2.641& 0.6674& 2.2605& 0.0105& 0.0289& 1.9344& 2.5989& 0.6645   \\
		&&$ \alpha_{1}$&
		3.7863& 0.2863& 0.6924& 2.2077& 5.3649& 3.1572& 3.6808& 0.1808& 0.6654& 2.3371& 5.4131& 3.076 \\
		&& $ \alpha_{2}$&3.7417& 0.2417& 0.5816& 2.2889& 5.1945& 2.9056& 3.6409& 0.1409& 0.5590& 2.3921& 5.2269& 2.8348  \\
		&&$ \alpha_{3}$& 3.7956& 0.2956& 0.6998& 2.2081& 5.3831& 3.175& 3.6923& 0.1923& 0.6724& 2.3403& 5.4355& 3.0952  \\
		&&$ \alpha_{4}$&3.7891& 0.2891& 0.6384& 2.2702& 5.308& 3.0378& 3.6870& 0.1870& 0.6147& 2.3855& 5.3476& 2.9621  \\
		&&$ \alpha_{5}$&3.7390& 0.2390& 0.5662& 2.3026& 5.1754& 2.8728& 3.6396& 0.1396& 0.5445& 2.4038& 5.2054& 2.8016  \\
		&&$ \alpha$& 3.5725& 0.0725& 0.2502& 2.6066& 4.5383& 1.9317& 3.4778& -0.0222& 0.2397& 2.6215& 4.5043& 1.8828  \\
		&& $R(t)$&0.5817& 0.0122& 0.0007& 0.5283& 0.6351& 0.1068& 0.5787& 0.0093& 0.0008& 0.5242& 0.6316& 0.1074  \\
		&& $h(t)$ &1.0287& -0.0187& 0.0066& 0.8699& 1.1874& 0.3175& 1.0205& -0.0268& 0.0066& 0.8675& 1.1842& 0.3167  \\
		& &$\text{MTF}$ &0.8984& 0.0248& 0.0029& 0.7926& 1.0043& 0.2117& 0.8972& 0.0235& 0.0031& 0.7930& 1.0106& 0.2176 \\
		225&3& $\beta$& 2.3088& 0.0588& 0.0291& 1.9750& 2.6426& 0.6676& 2.2622& 0.0122& 0.0289& 1.9363& 2.6003& 0.664\\
		&&$ \alpha_{1}$&3.7914& 0.2914& 0.6952& 2.2099& 5.3729& 3.1630& 3.6860& 0.1860& 0.6675& 2.3423& 5.4236& 3.0813   \\
		&&$ \alpha_{2}$&3.7838& 0.2838& 0.5947& 2.3125& 5.2551& 2.9427& 3.6829& 0.1829& 0.5709& 2.4184& 5.2857& 2.8673 \\
		&&$ \alpha_{3}$& 3.7913& 0.2913& 0.6979& 2.2061& 5.3764& 3.1703& 3.6881& 0.1881& 0.6698& 2.3369& 5.4251& 3.0882\\
		&&$ \alpha_{4}$& 3.7819& 0.2819& 0.6344& 2.2667& 5.2971& 3.0304& 3.6798& 0.1798& 0.6096& 2.3827& 5.3365& 2.9538\\
		&&$ \alpha_{5}$& 3.7541& 0.2541& 0.574& 2.3112& 5.197& 2.8858& 3.6542& 0.1542& 0.5518& 2.4118& 5.23& 2.8181\\
		&&$\alpha$&3.5848& 0.0848& 0.2526& 2.6148& 4.5548& 1.9400& 3.4900& -0.0100& 0.2416& 2.6306& 4.5212& 1.8906 \\
		&& $R(t)$&0.5811& 0.0117& 0.0007& 0.5277& 0.6345& 0.1068& 0.5782& 0.0088& 0.0008& 0.5236& 0.6311& 0.1075  \\
		&& $h(t)$ & 1.0311& -0.0162& 0.0066& 0.872& 1.1902& 0.3182& 1.0229& -0.0244& 0.0066& 0.8695& 1.1872& 0.3177  \\
		& &$\text{MTF}$ & 0.8971& 0.0234& 0.0029& 0.7915& 1.0026& 0.2111& 0.8958& 0.0222& 0.0031& 0.7919& 1.0091& 0.2172 \\
		\hline
	\end{tabular}
\end{table}

\begin{table}[H]
	\centering
	\small
	\caption{\small Point and interval estimates under the BAPCS setup 6.}
	\label{t7}
	\setlength\tabcolsep{1.5pt}
	\renewcommand{\arraystretch}{1.2}
	\begin{tabular}{cccccccccccccccccc}
		\hline
		&&&  \multicolumn{6}{c}{ML Estimates} & \multicolumn{6}{c}{~~ Pivotal Estimates} \\
		\cmidrule(lr){4-9}  \cmidrule(lr){10-17}
		&&&Est.&Bias& Variance &\multicolumn{2}{c}{ACIs}&L&    Est.&Bias& Variance &\multicolumn{2}{c}{~~ ~~GCIs~~}&L&&   \\
		\cmidrule(lr){7-8}  \cmidrule(lr){13-14}
		$n$&CP&Pars.&&&  &Lower & Upper  & & & &  & Lower &Upper\\ \hline
		250&1& $ \beta$& 2.2995& 0.0495& 0.0252& 1.9889& 2.6100& 0.6211& 2.2579& 0.0079& 0.0250& 1.9539& 2.5721& 0.6182   \\
		&&$ \alpha_{1}$&
		3.6951& 0.1951& 0.4769& 2.3737& 5.0165& 2.6429& 3.6083& 0.1083& 0.4605& 2.4561& 5.0396& 2.5835 \\
		&& $ \alpha_{2}$&3.7036& 0.2036& 0.4406& 2.4311& 4.9761& 2.545& 3.6172& 0.1172& 0.4256& 2.5024& 4.9949& 2.4925  \\
		&&$ \alpha_{3}$& 3.7156& 0.2156& 0.5430& 2.3074& 5.1237& 2.8163& 3.6268& 0.1268& 0.5246& 2.4083& 5.1638& 2.7555 \\
		&&$ \alpha_{4}$&3.7463& 0.2463& 0.5568& 2.3250& 5.1675& 2.8424& 3.6554& 0.1554& 0.5377& 2.4276& 5.2035& 2.7758 \\
		&&$ \alpha_{5}$&3.7390& 0.2390& 0.5662& 2.3026& 5.1754& 2.8728& 3.6396& 0.1396& 0.5445& 2.4038& 5.2054& 2.8016  \\
		&&$ \alpha$& 3.5577& 0.0577& 0.2032& 2.6859& 4.4295& 1.7436& 3.475& -0.025& 0.1952& 2.6943& 4.398& 1.7037  \\
		&& $R(t)$&0.5806& 0.0112& 0.0006& 0.5312& 0.63& 0.0988& 0.5777& 0.0082& 0.0006& 0.5273& 0.6266& 0.0993  \\
		&& $h(t)$ &1.0293& -0.0180& 0.0053& 0.887& 1.1716& 0.2846& 1.0228& -0.0245& 0.0053& 0.885& 1.1691& 0.2841  \\
		& &$\text{MTF}$ &0.8960& 0.0224& 0.0025& 0.7991& 0.9929& 0.1939& 0.894& 0.0203& 0.0026& 0.7983& 0.9967& 0.1983 \\
		250&2& $\beta$& 2.2958& 0.0458& 0.024& 1.9927& 2.599& 0.6063& 2.2578& 0.0078& 0.0239& 1.9606& 2.5649& 0.6043& \\
		&&$ \alpha_{1}$&3.7003& 0.2003& 0.462& 2.3967& 5.0039& 2.6072& 3.6226& 0.1226& 0.4484& 2.4771& 5.0336& 2.5565&  \\
		&&$ \alpha_{2}$&3.6819& 0.1819& 0.4201& 2.4371& 4.9267& 2.4896& 3.6044& 0.1044& 0.4074& 2.5088& 4.9484& 2.4396& \\
		&&$ \alpha_{3}$& 3.7247& 0.2247& 0.5276& 2.3338& 5.1155& 2.7817& 3.645& 0.145& 0.512& 2.4323& 5.1589& 2.7266& \\
		&&$ \alpha_{4}$& 3.7129& 0.2129& 0.4668& 2.4039& 5.022& 2.6181& 3.6345& 0.1345& 0.4524& 2.4873& 5.0523& 2.5651& \\
		&&$ \alpha_{5}$& 3.7334& 0.2334& 0.5314& 2.3411& 5.1256& 2.7845& 3.6517& 0.1517& 0.5148& 2.4389& 5.1669& 2.728& \\
		&&$\alpha$&3.5568& 0.0568& 0.1896& 2.7134& 4.4003& 1.6869& 3.4822& -0.0178& 0.1832& 2.7202& 4.373& 1.6528& \\
		&& $R(t)$&0.5793& 0.0098& 0.0006& 0.5294& 0.6291& 0.0997& 0.5765& 0.007& 0.0007& 0.5256& 0.626& 0.1004&   \\
		&& $h(t)$ & 1.032& -0.0154& 0.0052& 0.8905& 1.1734& 0.2828& 1.0263& -0.0210& 0.0053& 0.8892& 1.1718& 0.2826&   \\
		& &$\text{MTF}$ & 0.8934& 0.0198& 0.0025& 0.7959& 0.9909& 0.1950& 0.8915& 0.0178& 0.0026& 0.7952& 0.9947& 0.1995&  \\
		250&3& $\beta$& 2.2929& 0.0429& 0.0240& 1.9900& 2.5958& 0.6058& 2.2550& 0.0050& 0.0238& 1.9582& 2.5615& 0.6033& \\
		&&$ \alpha_{1}$&3.6977& 0.1977& 0.4635& 2.3945& 5.001& 2.6066& 3.6199& 0.1199& 0.4486& 2.4774& 5.0291& 2.5516& 
		&  \\
		&&$ \alpha_{2}$&3.7020& 0.2020& 0.4264& 2.4492& 4.9549& 2.5057& 3.6243& 0.1243& 0.4132& 2.5225& 4.9795& 2.4570&  \\
		&&$ \alpha_{3}$& 3.7186& 0.2186& 0.5290& 2.3297& 5.1076& 2.7779& 3.6384& 0.1384& 0.5127& 2.4297& 5.1479& 2.7182&  \\
		&&$ \alpha_{4}$& 3.6814& 0.1814& 0.4584& 2.3847& 4.9780& 2.5933& 3.6044& 0.1044& 0.4437& 2.4678& 5.0095& 2.5416&  \\
		&&$ \alpha_{5}$& 3.7198& 0.2198& 0.5250& 2.3339& 5.1058& 2.7720& 3.6390& 0.1390& 0.5083& 2.4314& 5.1459& 2.7145&  \\
		&&$\alpha$&3.5512& 0.0512& 0.1893& 2.7091& 4.3933& 1.6843& 3.4771& -0.0229& 0.1825& 2.7172& 4.366& 1.6488& \\
		&& $R(t)$& 0.5791& 0.0097& 0.0006& 0.5293& 0.629& 0.0997& 0.5763& 0.0069& 0.0007& 0.5255& 0.6257& 0.1002&  \\
		&& $h(t)$ & 1.0317& -0.0156& 0.0052& 0.8903& 1.1731& 0.2828& 1.0261& -0.0212& 0.0052& 0.8890& 1.1717& 0.2827&   \\
		& &$\text{MTF}$ &  0.8934& 0.0197& 0.0025& 0.7958& 0.991& 0.1952& 0.8914& 0.0177& 0.0026& 0.7951& 0.9945& 0.1994& \\
		
		\hline
	\end{tabular}
\end{table}

\newpage

\section{Real Data Analysis} \label{sec 7}

In this section, we analyze the flexibility of the IEP distribution by comparing it with other four lifetime distributions, namely, generalized Pareto (GP) distribution, exponentiated Pareto (EP) distribution, inverted exponentiated Rayleigh (IER) distribution and inverse Lomax (IL) distribution. The PDFs of the above four proposed distributions are as follows:

 \begin{itemize}
 	\item GP: \begin{equation}\nonumber
 		f_{GP}(x; k, \sigma) = \begin{cases}
 			\frac{1}{\sigma}(1 - \frac{kx}{\sigma})^{\frac{1}{k} - 1}; & \text{if $x>0$, $k \neq 0$, $\sigma > 0,$}\\
 		\frac{1}{\sigma}\exp(-x/\sigma); & \text{if $x>0$, $k = 0$, $\sigma > 0$}.
 		\end{cases}
 	\end{equation}
 	\item EP: \begin{equation}\nonumber
 	  f_{EP}(x; \lambda, \theta) = \lambda \theta \left[ 1 - (1 + x)^{-\lambda} \right]^{\theta - 1}(1+x)^{-(\lambda+1)};~x>0,~\lambda, \theta > 0.
  \end{equation}
   \item IER: \begin{equation} \nonumber
   	f_{IER}(x; \alpha, \beta) = 2\alpha \beta x^{-3}\exp\left(-\frac{\beta}{x^2}\right)\left[1 - \exp\left(-\frac{\beta}{x^2}\right)\right]^{\alpha-1};~x>0,~\alpha, \beta > 0.
   \end{equation}
    
    \item IL: \begin{equation}\nonumber
    \displaystyle{
    f_{IL}(x; \alpha, \theta) = \frac{\alpha}{\theta}\frac{1}{x^2}\frac{1}{(1 + \frac{1}{\theta x})^{\alpha + 1}}};~x>0,~\alpha, \theta > 0.
    \end{equation}
 \end{itemize}

A real dataset consisting of strengths measured in GPa, of 69 single carbon fibres tested under tension at gauge lengths of 20 mm, is taken from \cite{ahmad2021transmuted}. The strengths are given as follows:\\

\noindent 0.312, 0.314, 0.479, 0.552, 0.700, 0.803, 0.861, 0.865, 0.944, 0.958, 0.966, 0.977, 1.006, 1.021, 1.027, 1.055, 1.063, 1.098, 1.140, 1.179, 1.224, 1.240, 1.253, 1.270, 1.272, 1.274, 1.301, 1.301, 1.359, 1.382, 1.382, 1.426, 1.434, 1.435, 1.478, 1.490, 1.511, 1.514, 1.535, 1.554, 1.566, 1.570, 1.586, 1.629, 1.633, 1.642, 1.648, 1.684, 1.697, 1.726, 1.770, 1.773, 1.800, 1.809, 1.818, 1.821, 1.848, 1.880, 1.954, 2.012, 2.067, 2.084, 2.090, 2.096, 2.128, 2.233, 2.433, 2.585, 2.585.

\begin{figure}[H]
   	\centering
   	{\includegraphics[height = 11cm, width = 11cm]{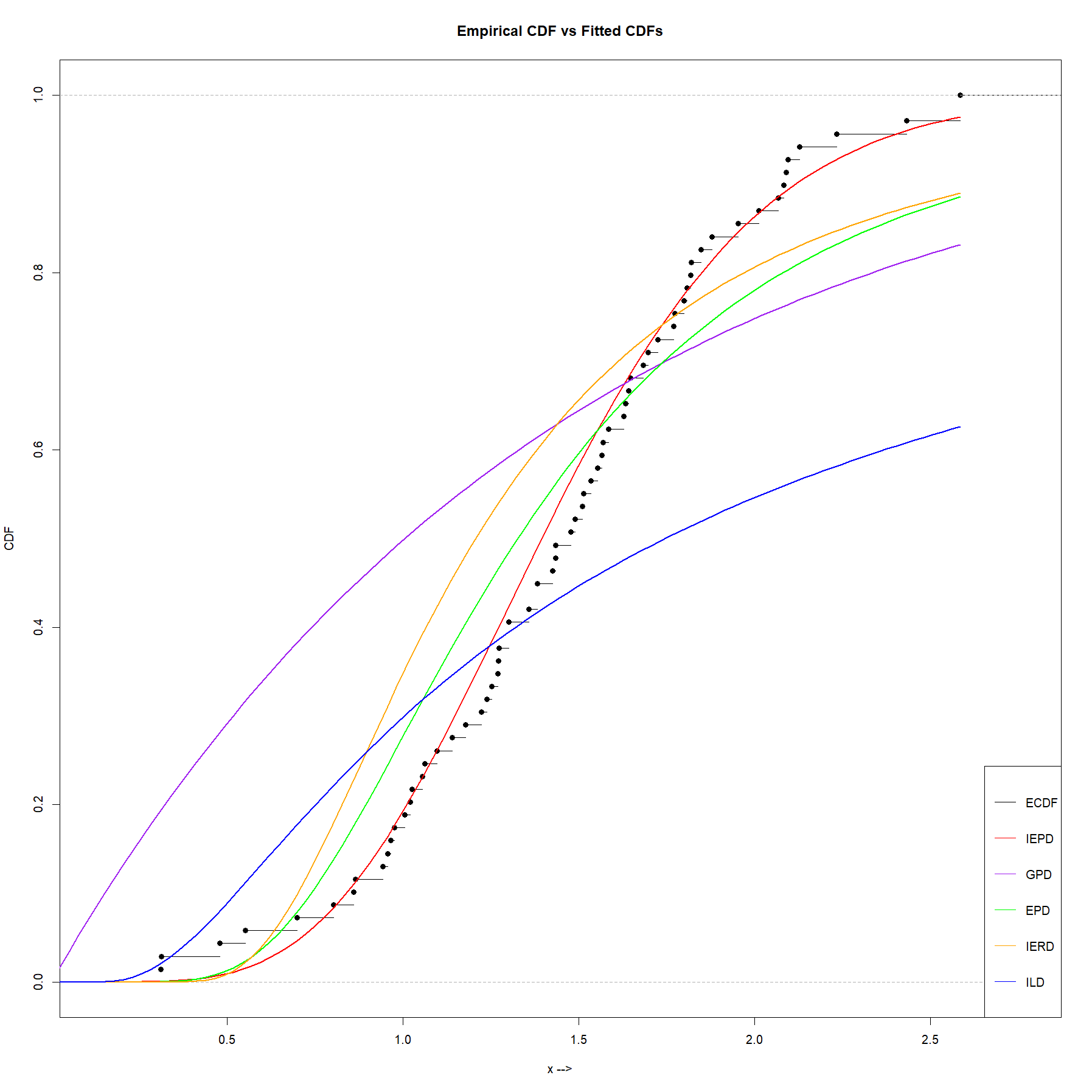}}
   	\caption{Plot of the empirical CDF of the real data along with the proposed theoretical CDFs.}
   	\label{ecdf}
   \end{figure}

\begin{figure}[H]
   	\centering
   	{\includegraphics[height = 11cm, width = 11cm]{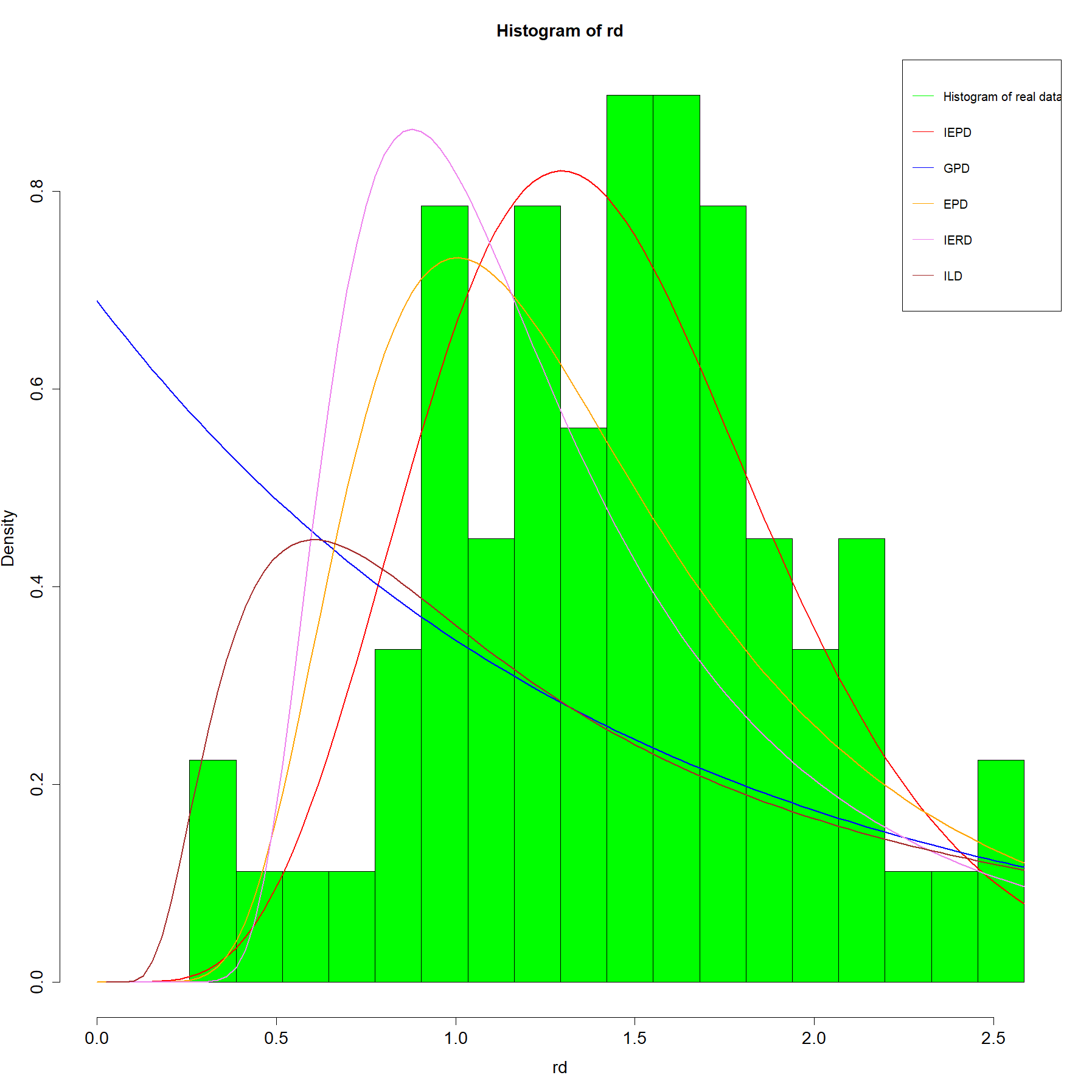}}
   	\caption{Plot of the histogram of the real data along with the proposed theoretical PDFs.}
   	\label{hist}
   \end{figure}

\begin{figure}[H]
	\centering
	\subfigure[]{ \includegraphics[height = 6cm, width = 6cm]{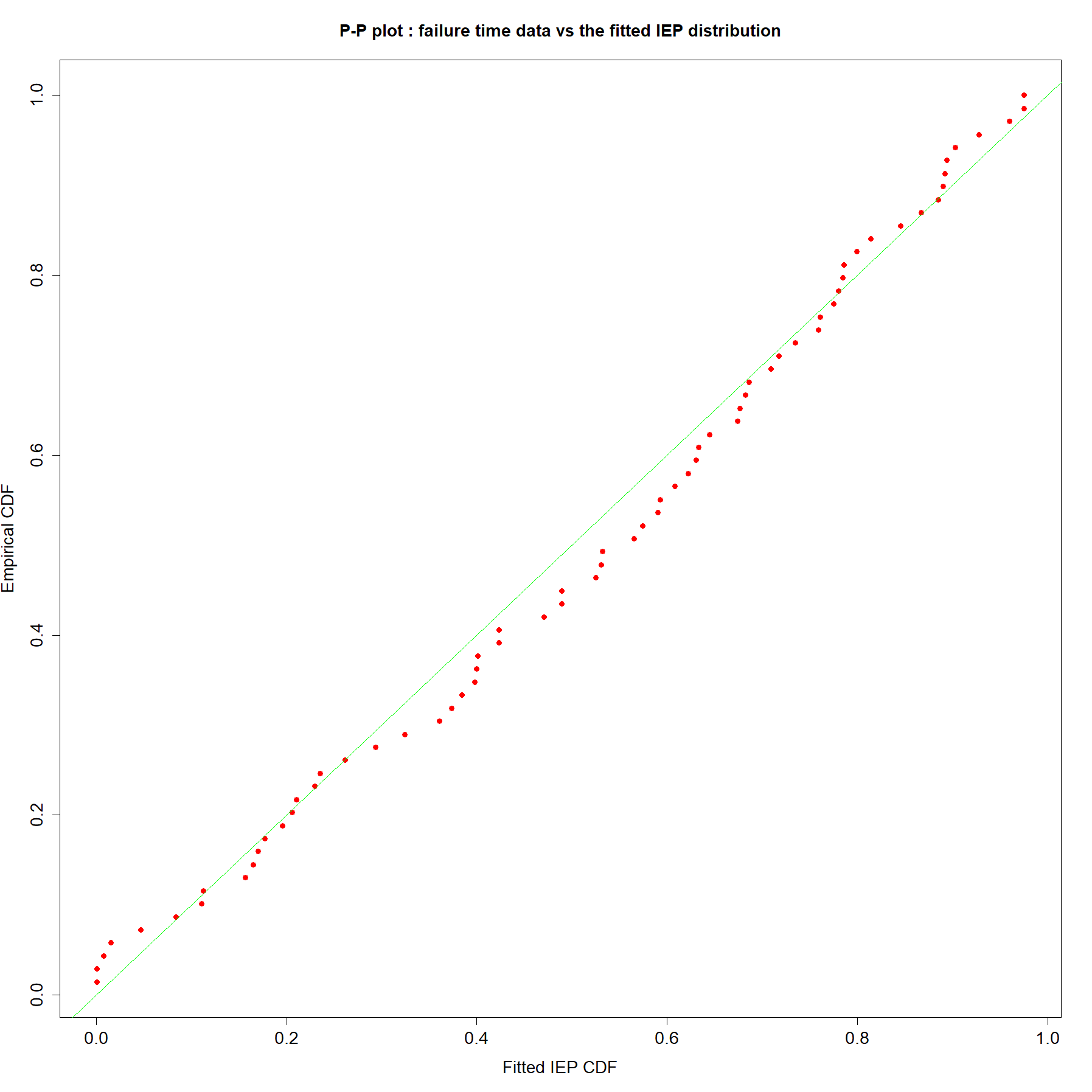}}
	\subfigure[]{\includegraphics[height = 6cm, width = 6cm]{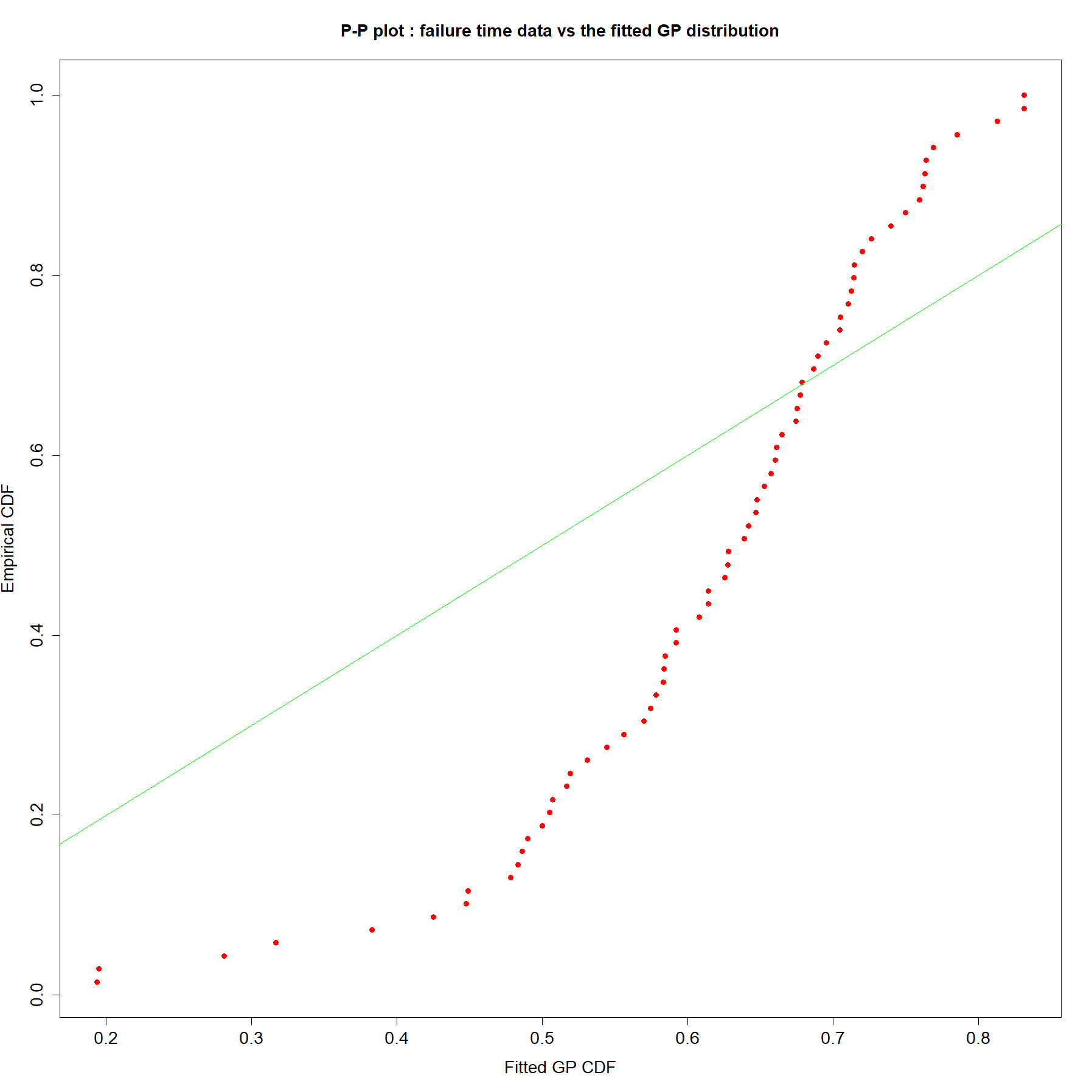}}
	\subfigure[]{\includegraphics[height = 6cm, width = 6cm]{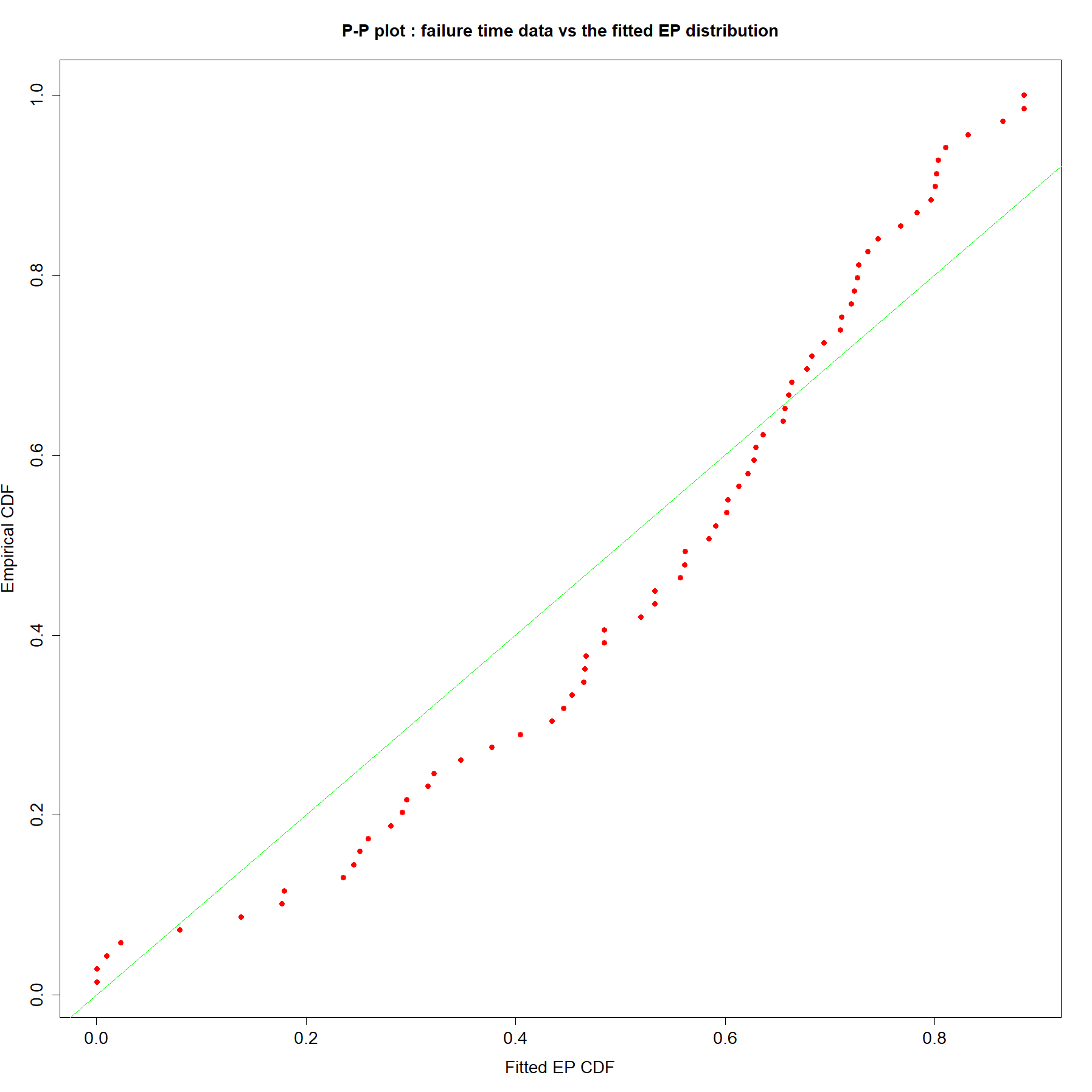}}
	\subfigure[]{ \includegraphics[height = 6cm, width = 6cm]{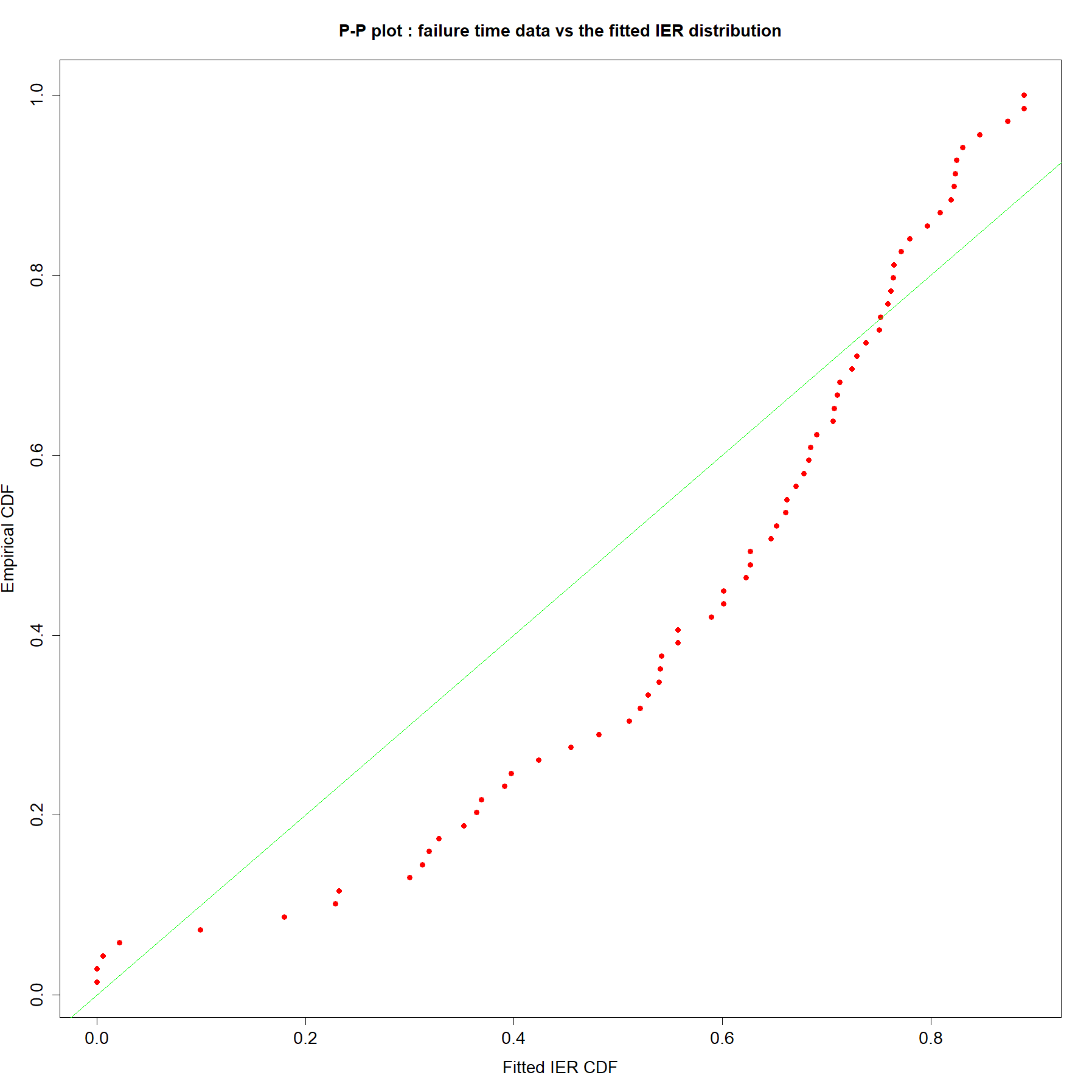}}
	\subfigure[]{ \includegraphics[height = 6cm, width = 6cm]{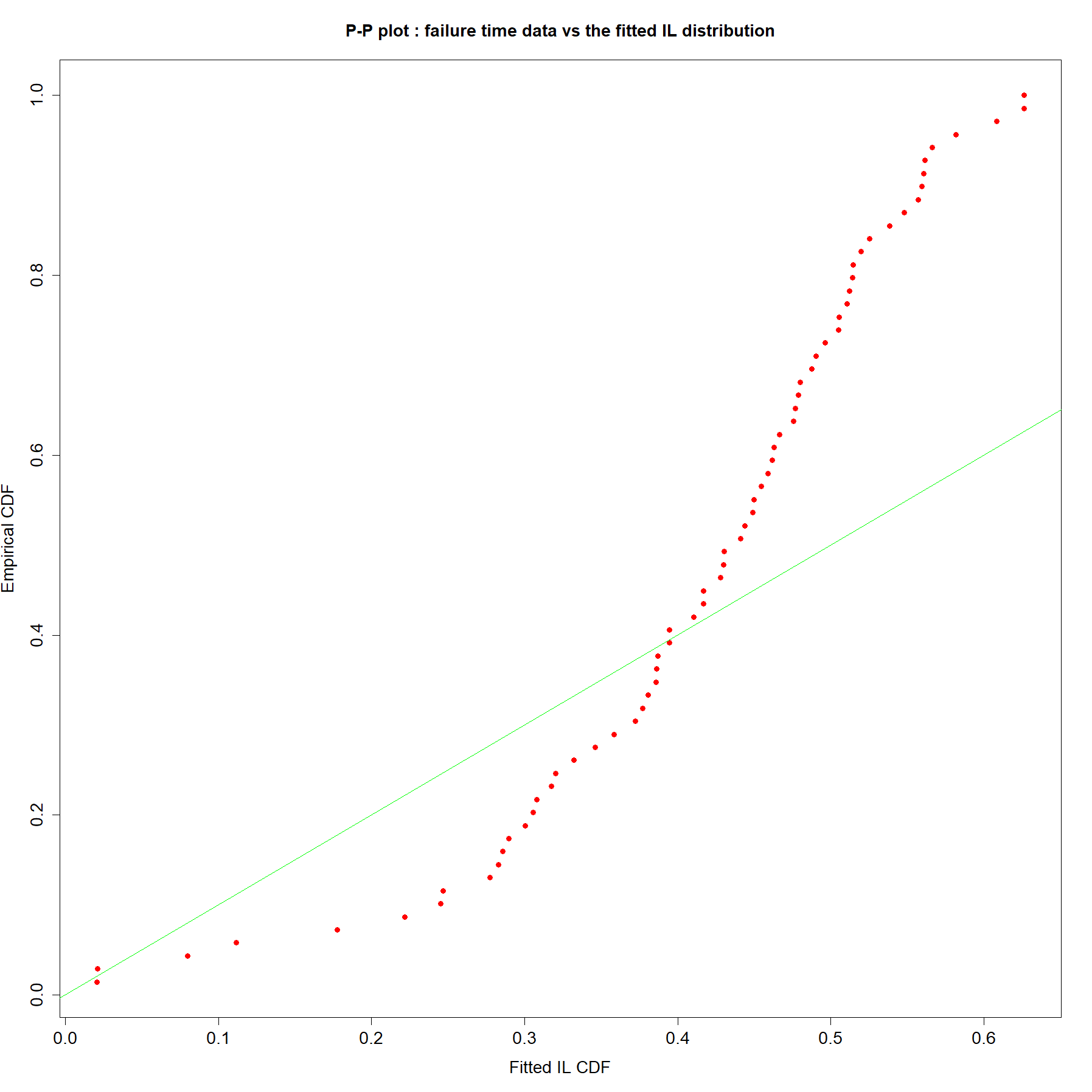}}
	\caption{P-P plot of the dataset vs fitted (a) IEP distribution, (b) GP distribution, (c) EP distribution, (d) IER distribution, and (e) IL distribution.}
	\label{fig 2}	
\end{figure}

\begin{figure}[H]
	\centering
	\subfigure[]{ \includegraphics[height = 6cm, width = 6cm]{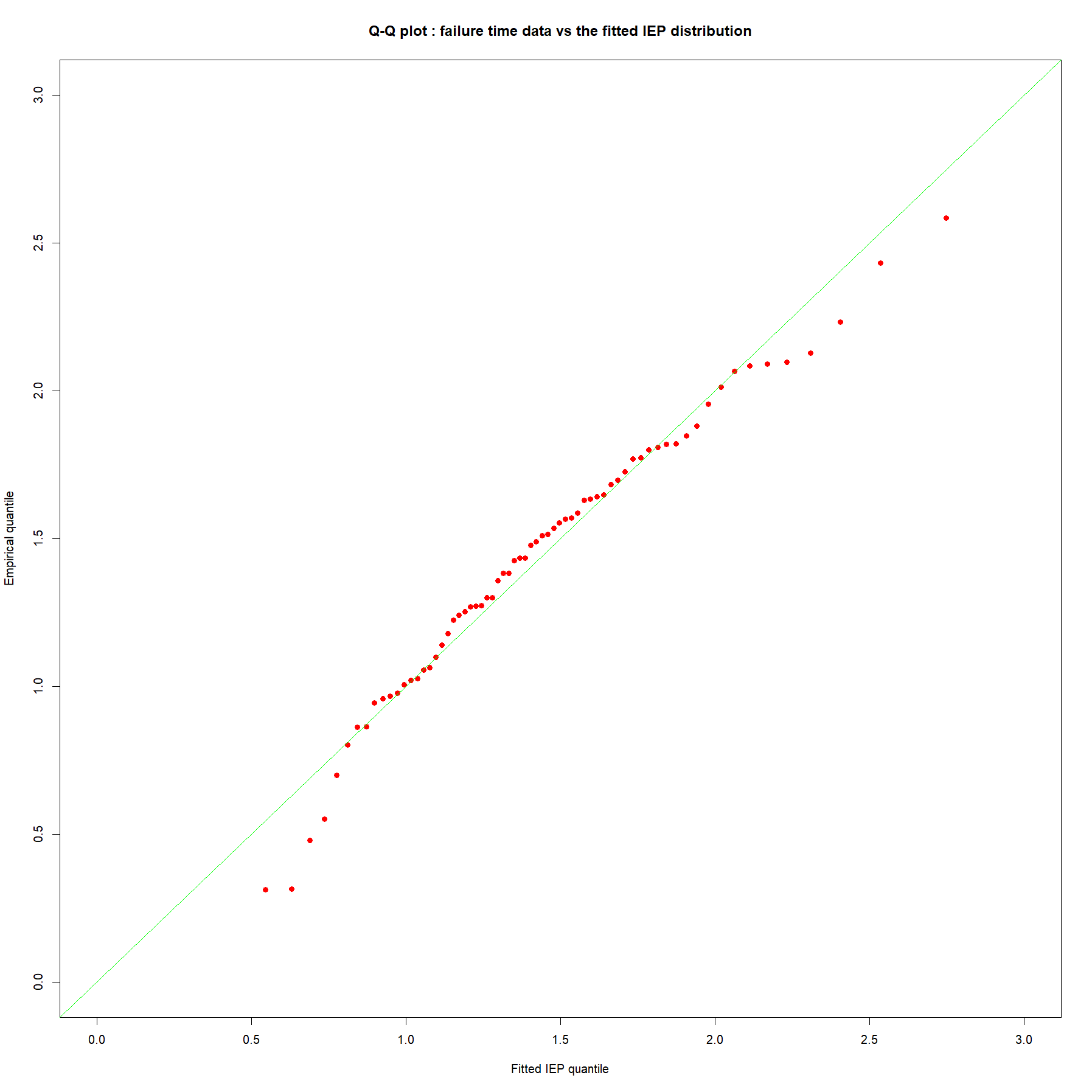}}
	\subfigure[]{\includegraphics[height = 6cm, width = 6cm]{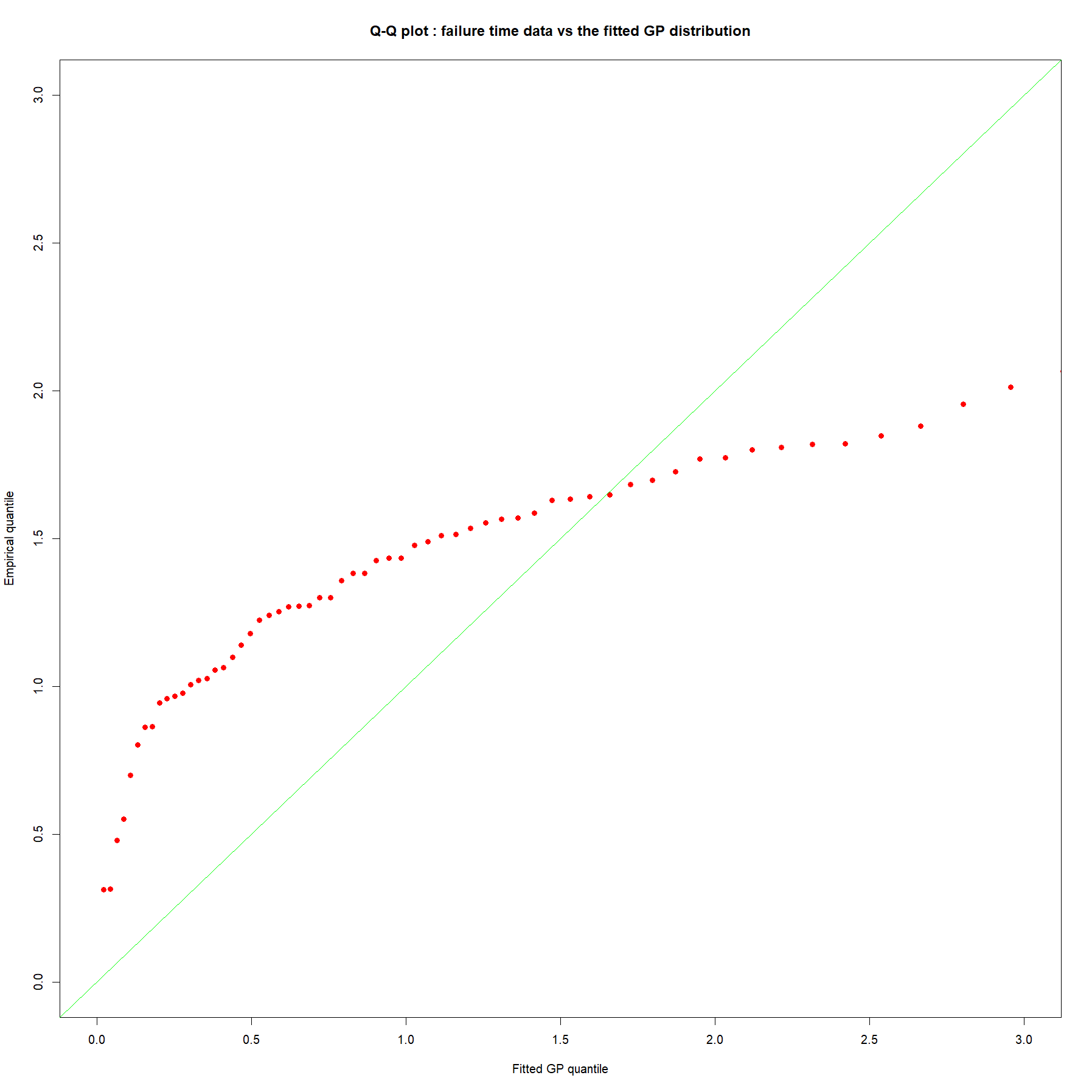}}
	\subfigure[]{	\label{figure07} \includegraphics[height = 6cm, width = 6cm]{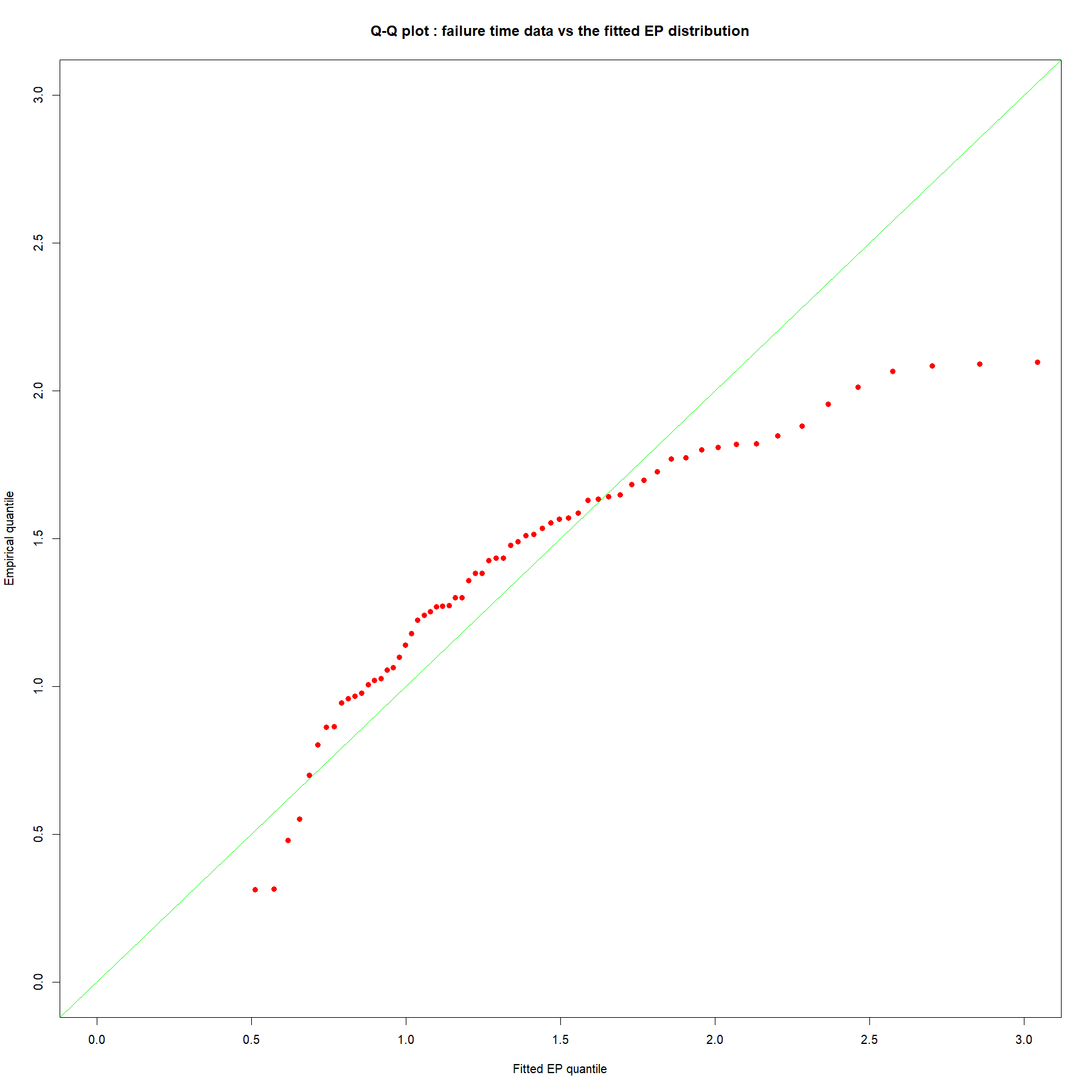}}
	\subfigure[]{\includegraphics[height = 6cm, width = 6cm]{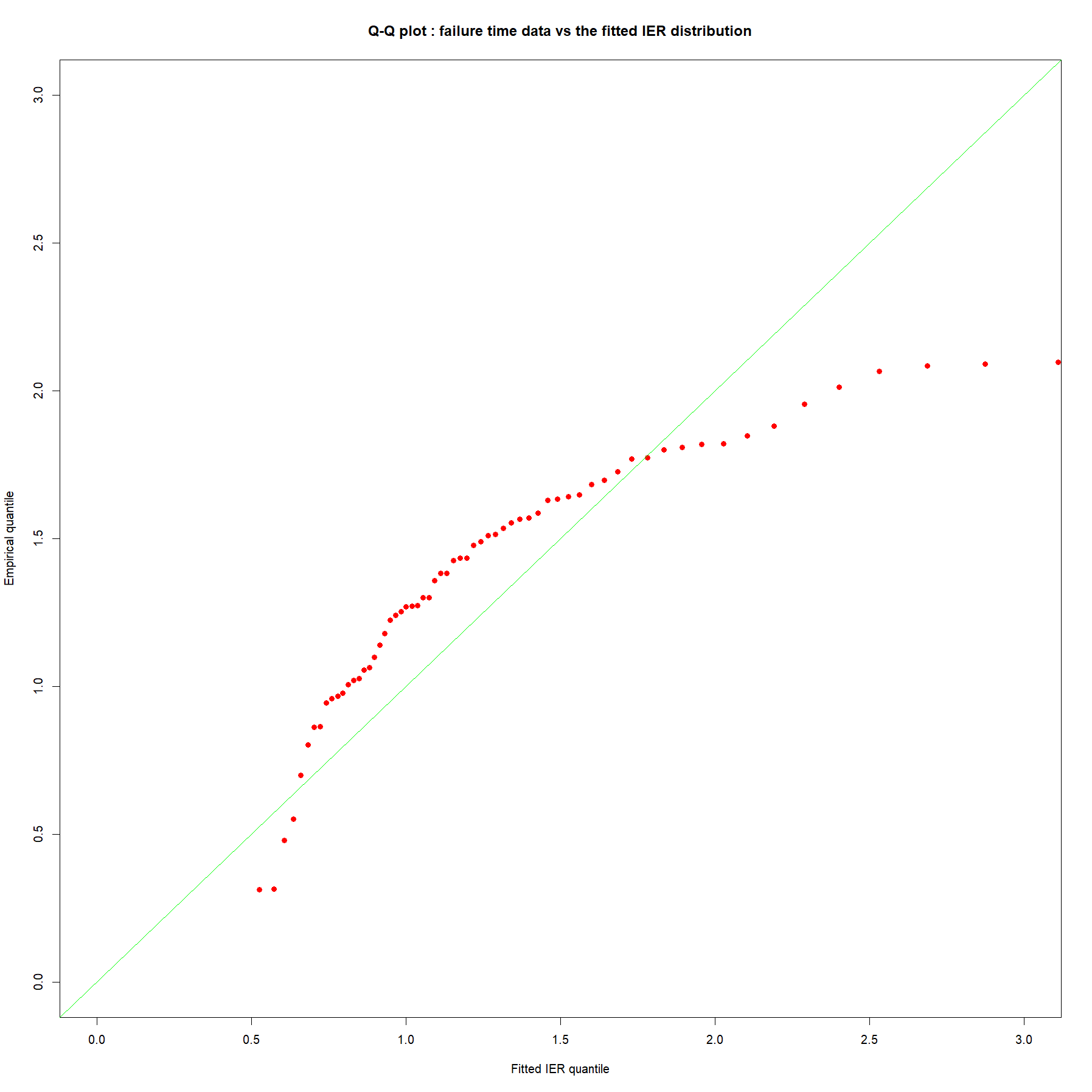}}
	\subfigure[]{\label{figure08} \includegraphics[height = 6cm, width = 6cm]{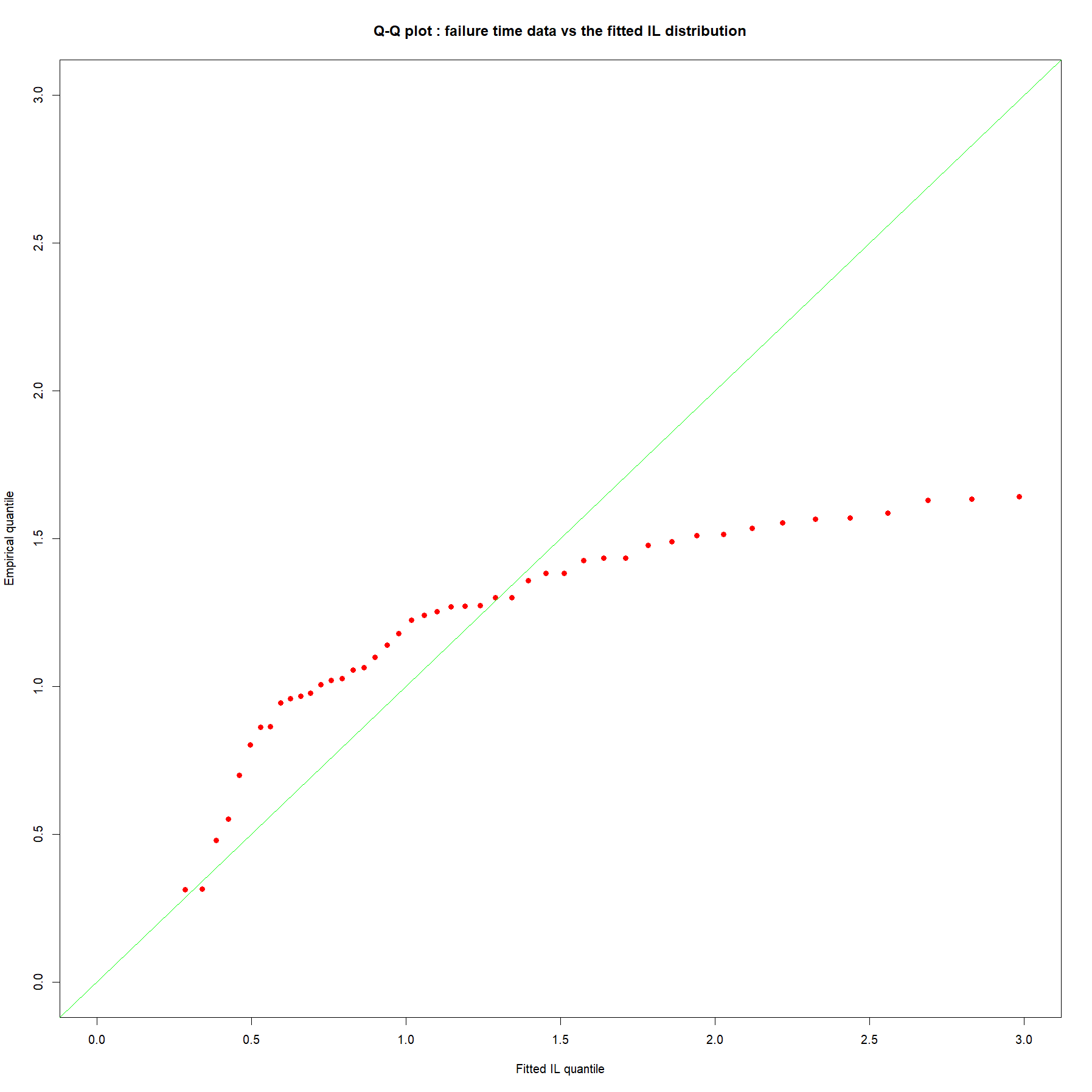}}
	\caption{Q-Q plot of the dataset vs fitted (a) IEP distribution, (b) GP distribution, (c) EP distribution, (d) IER distribution, and (e) IL distribution.}
	\label{fig 3}	
\end{figure}

\begin{figure}[H]
	\centering
	{ \includegraphics[height = 12cm, width = 12cm]{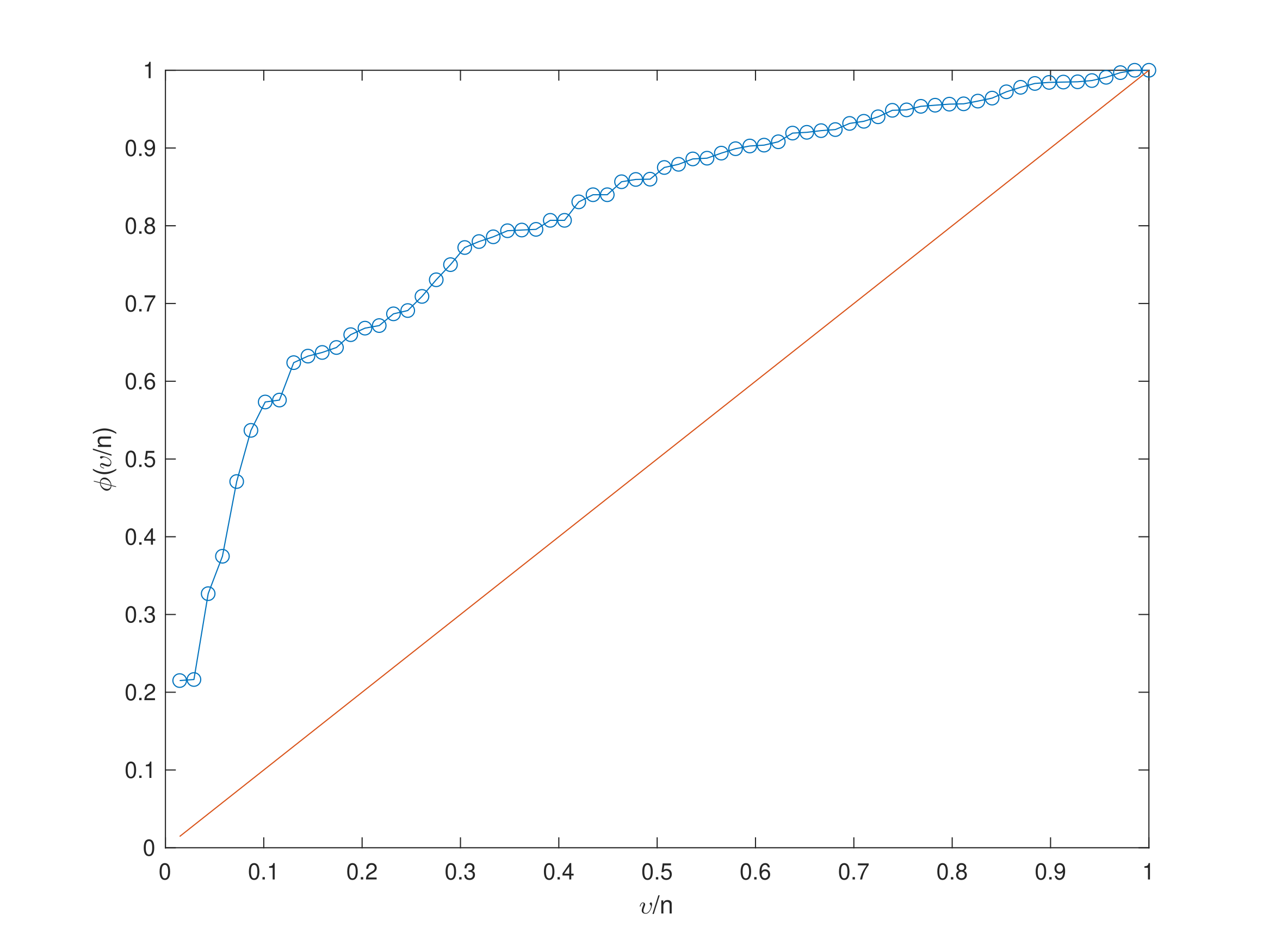}}
	\caption{Plot of the scaled total time on test based on the real data.}
	\label{TTT}
\end{figure}

 The unknown parameters are estimated by the MLE method. The estimated parameters, values of the Kolmogorov-Smirnov statistics (K-S statistics), $p$-value in the K-S test, values of the Akaike information criterion (AIC), Bayesian information criterion (BIC), consistent Akaike information criterion (CAIC) and Hannan-Quin information criterion (HQIC) are listed in Table \ref{Table 1}. All the information criteria and the value of the K-S statistic is smallest when the IEP distribution is in use. Also, the $p$-value of the K-S test is largest in case of the IEP model.  The plot of the empirical CDF along with the CDF plots for all the proposed distributions are given in Figure \ref{ecdf}. Moreover, the histogram plot of the real data along with the theoretical PDFs is displayed in the Figure \ref{hist}. The P-P plots and the Q-Q plots are displayed by the figure \ref{fig 2} and \ref{fig 3} respectively. It is evident from all these plots and the Table \ref{Table 1} that the IEP distribution model is more appropriate than the other four models to fit the carbon fibre strength data. Figure \ref{TTT} is the total time on test plot which have a concave downward shape, suggests that the hazard rate of the real data mentioned in this section is increasing.

\begin{table}[H]
	\centering
	\caption{\label{Table 1} MLEs of the parameters using the carbon fibre test data and AIC, BIC, CAIC, HQIC measures, {K-S} statistics and $p$-value.}
	\label{t23}
	\setlength\tabcolsep{1.5pt}
	\scalebox{0.945}{
	\begin{tabular}{|c|c|c|c|c|c|c|c|c|}
		\hline
		Model & Pars. & MLE & AIC & BIC & CAIC & HQIC & K-S statistic & $p$-value\\
		\hline
		IEP &$(\alpha, \beta)$  & (43.8478, 7.6876) &  108.6455 & 113.1138 & 115.1138 & 110.4182 & 0.0755 & 0.8266\\
		
		GP &$(k, \sigma)$  & (0, 1.4510) & 193.3750 & 197.8433 & 199.8433 &  195.1477 &  0.3623 & 2.716 $\times 10^{-8}$ \\
		
		EP & $(\lambda, \theta)$  & (3.9902, 19.7892) & 134.6356 & 139.1038 &141.1038 &  136.4083 & 0.1451 & 0.1095\\
		
		IER &$(\alpha, \beta)$  & (1.2358, 1.2322) &  160.1693 & 164.6375 & 166.6375 & 161.9420 & 0.2210 &0.0024\\
		
		IL &$(\alpha, \theta)$  & (52679.34,  43536.79) & 198.2299 &  202.6981 & 204.6981 & 200.0026 & 0.3757 & $6.933 \times 10^{-9}$\\
		\hline
		
	\end{tabular}}
\end{table}

\section{Concluding Remarks} \label{sec 8}
This article delves into the parameter estimation of the IEP distribution and estimation of some of its RCs such as, the reliability or suvival function, the hazard rate function, and the median time to failure. The block censoring mechanism is considered to generate the samples in the simulation study. The adaptive progressive Type-II censoring scheme is combined with the block censoring procedure which is a new idea in the study of block progressive censoring. In the simulation study, the MLE method and the pivotal estimation method were used to carry out both the point and interval estimations of the unknown parameters and the proposed RCs of the IEP distribution. As interval estimation, the asymptotic confidence intervals and the generalized confidence intervals are computed. In addition, the idea of order statistics is also used to explore various mathematical developments of the IEP distribution. All the estimated values and estimated intervals are tabulated for thorough obsevations. It was inetrestingly revealed that the pivotal method gives better eatimations than the classical MLE method and this fact can give motivation to implement the pivotal estimation method in the study of staistical inference with block censoring. In the section of real data analysis, the IEP distribution is compared with other four well-known distributions by considering a carbon fibre strength dataset. Various statistical plots are presented in this section to compare among the goodnesses of fit of the  proposed lifetime distributions with the real data and it was noticed that the IEP distribution has the better fit compared to the other four proposed distributions. Moreover, various information criteria are computed corresponding to all the five proposed distributions and the K-S test is also done to check the goodnesses of fit and then it is again evident that the IEP distribution is better fitted distribution than each of the four others.

\section*{Acknowledgements}

The financial support (institute fellowship, Indian Institute of Technology Roorkee, Govt. of India), is sincerely acknowledged with thanks by Rajendranath Mondal. Dr. Raju Bhakta gratefully acknowledges the financial support for this research work under the project grant IITR/SRIC/2301/DRD-2236-CSE/23-24/SR-03.

\section*{Conflict of Interest} All the authors have no conflict of interest.

 \bibliography{Conference_paper}  
\end{document}